\documentclass[12pt,draftcls,onecolumn]{IEEEtran}
\usepackage{graphicx}
\usepackage{amstext}
\usepackage{amsmath}
\usepackage{amsfonts,amssymb}
\usepackage{amssymb}
\usepackage[dvips]{epsfig}
\usepackage{epsfig}
\usepackage[dvips]{graphicx}
\usepackage{multirow}
\usepackage{multirow} 
\usepackage[utf8]{inputenc}
\usepackage[english]{babel}
\usepackage{caption}
\usepackage{subcaption}
\usepackage{float}
\usepackage{ltablex}
\usepackage{algorithmic}
\usepackage{setspace}
\newtheorem{theorem}{Theorem}[section]

\newtheorem{lemma}[theorem]{Lemma}

\def\bv{{\mathbf{V}}}
\def\ba{{\mathbf{A}}}

\def\bvv{{\mathbf{v}}}
\def\bva{{\mathbf{v}}}
\def\bw{{\mathbf{W}}}

\def\bu{{\mathbf{u}}}
\def\bX{{\mathbf{X}}}
\def\bwv{{\mathbf{w}}}
\def\bhv{{\mathbf{h}}}
\def\bh{{\mathbf{H}}}

\def\bx{{\mathbf{X}}}
\def\bxv{{\mathbf{x}}}

\def\bz{{\mathbf{Z}}}

\providecommand{\keywords}[1]{\textbf{\textit{Index terms---}} #1}
\begin{document}

\title{Novel Algorithms based on Majorization Minimization for Nonnegative Matrix Factorization}
\author{R.~Jyothi, P.~Babu and R.~Bahl
\footnote{The authors are with CARE, IIT Delhi, New Delhi, 110016, India.(email: jyothi.r@care.iitd.ac.in, prabhubabu@care.iitd.ac.in, rbahl@care.iitd.ac.in)}}
\maketitle
\maketitle
\begin{abstract}
Matrix decomposition is ubiquitous and has applications in various fields like speech processing, data mining and image processing to name a few. Under matrix decomposition, nonnegative matrix factorization is used to decompose a nonnegative matrix into a product of two nonnegative matrices which gives some meaningful interpretation of the data. Thus, nonnegative matrix factorization has an edge over the other decomposition techniques. In this paper, we propose two novel iterative algorithms based on Majorization Minimization (MM)-in which we formulate a novel upper bound and minimize it to get a closed form solution at every iteration. Since the algorithms are based on MM, it is ensured that the proposed methods will be monotonic. The proposed algorithms differ in the updating approach of the two nonnegative matrices. The first  algorithm-\textbf{I}terative \textbf{No}nnegative \textbf{M}atrix Factorization (\textbf{INOM}) sequentially updates the two nonnegative matrices while the second algorithm-\textbf{Par}allel \textbf{I}terative \textbf{No}nnegative \textbf{M}atrix Factorization (\textbf{PARINOM}) parallely updates them. We also prove that the proposed algorithms converge to the stationary point of the problem. Simulations were conducted to compare the proposed methods with the existing ones and was found that the proposed algorithms performs better than the existing ones in terms of computational speed and convergence. 
\end{abstract}
\section{Introduction} \label{sec:1}
Recent advancements in sensor technology and communications has created huge collection of data resulting in big data matrices, which when appropriately analyzed can give useful insights about the data. Many times, researchers reduce the dimension of the data matrix for easier visualization and to lessen the computational load \cite{app_reduction}. There are various tools to reduce the dimension of the data, some of them are - \emph{Singular Value Decomposition (SVD)} \cite{svd}, \emph{Principal Component Analysis (PCA)} \cite{pca} and \emph{Factor Analysis (FA)} \cite{fa}. However, the major shortcoming of these techniques is that when the input data matrix is nonnegative as in the case of speech/image processing, the reduced data matrix can have negative entries, which makes the interpretation difficult. 
\emph{Nonnegative matrix factorization (NMF)} is another dimension reduction technique, which as the name suggests decomposes the matrix such that the reduced dimension matrices are always nonnegative. NMF has found applications in music processing \cite{music}, data mining \cite{text},  image processing \cite{image} and in neurobiology \cite{brain} to name a few. Mathematically, NMF problem can be written as:
\begin{equation} \label{eq:11}
\begin{array}{ll}
\textrm{NMF:}\quad  \underset{\bw,\bh \geq 0}{\rm minimize}\:\{f_{_{\rm NMF}}\left(\bw,\bh\right)\overset{\Delta} = \: \|\bv -\bw\bh\|_F^2\}
\end{array}
\end{equation}
where, ${\| \bx \|_{F}}$ denotes the Frobenious norm of matrix ${\bx}$, ${\bv}$  is a data matrix made of real entries and of size ${n \times m}$, ${\bw}$ and ${\bh}$ are decomposed matrices of size $n \times r$  and $r \times m$, respectively. Here, $r$ is chosen to be lesser than $m$ and $n$. The constraint ($\bw,\bh \geq 0$) is such that $\bw$ and $\bh$ matrices must not have any nonnegative element.  The $i^{th}$ column of $\bv$ matrix, represented as $\bva_{i}$, can be expressed as nonnegative weighted linear combination of columns of $\bw$ matrix:
\begin{equation}\label{cone}
\bva_{i}= \displaystyle\sum_{j=1}^{r}h_{ji}\bwv_{j} 
\end{equation} 
where $\bwv_{j}$ is the  ${j^{th}}$ column of $\bw$ matrix and $h_{ji}$ is the $(j, i)^{th}$ element of $\bh$ matrix. Geometrically, this means that the columns of $\bw$ matrix generates a \emph{simplicial convex cone} (\cite{geo}, \cite{cone}) that contains the data matrix $\bva_{i}$, where \emph{simplicial convex cone C}  is defined as:
\begin{equation}
C = \left\{\displaystyle\sum_{i=1}^{r}\theta_{i}\bu_{i}:\theta_{i}\geq0 \right\}
\end{equation}
where $\theta_{i} \in \mathbf{R^{+}}$ and $\bu_{i}$ is a vector in $\mathbf{R}$ of arbitrary dimension. From (\ref{cone}), it can be seen that there exists many cones (whose vertices are defined by columns of $\bw$), containing the data $\bva_{i}$. In order to have a unique cone, one must normalize the columns of $\bw$ and $\bv$ matrix \cite{cone} and hence in this paper we do the same. 
The NMF problem in (\ref{eq:11}) is not jointly convex in ${\bw}$ and ${\bh}$. However, it is separately convex in either ${\bw}$ or ${\bh}$. Hence, alternatingly minimizing over ${\bh}$ and ${\bw}$ seems to be a favorable direction to solve the NMF problem in (\ref {eq:11}). This approach is traditionally named as \emph{Alternating minimization}. In alternating minimization, constrained minimization is performed with respect to one matrix while keeping the other matrix fixed and vice-versa. The pseudo code of Alternating minimization is shown below:
\begin{center}
\begin{tabular}{   l }
\hline
\hline
{\bf{Table 1: Pseudocode of Alternating Minimization}} (\cite{bcd1}, \cite{bcd2} and \cite{bcd3}) \\
\hline
\hline
{\bf{Input}}: Data sample $\bv$, with each column normalized, \emph{r}, \emph{m} and \emph{n}. \\
{\bf{Initialize}}: Set \emph{k} = 0. Initialize ${\bw^{0}}$ and ${\bh^{0}}$. Each column of ${\bw^{0}}$ normalized. \\
{\bf{Repeat}}: \\
     1) Fix ${\bw}^{k}$ and find ${\bh^{k+1}}$= $\underset{\bh \geq 0}{\rm arg\,minimize} \: f_{_{\rm NMF}}\left(\bw^{k},\bh\right)$\\
     2) Fix ${\bh}^{k+1}$ and find ${\bw^{k+1}}$ = $\underset{\bw \geq 0}{\rm arg\,minimize} \:f_{_{\rm NMF}}\left(\bw,\bh^{k+1}\right)$\\
     3) Normalize the columns of ${\bw}^{k+1}$\\
    $k \leftarrow k+1$\\
{\bf{until convergence}}\\
\hline
\hline
\end{tabular}
\end{center}
Many algorithms have been proposed to solve the NMF problem. Some of them are discussed as follows: The baseline algorithm used to solve the above problem is the \emph{Multiplicative Update(MU)} algorithm proposed by Lee et. al. \cite{MU}. This algorithm which is iterative in nature, is  based on \emph{Block Majorization Minimization (MM)} (which we will introduce shortly in section \ref{sec:2}). The final update equations of both the matrices has only multiplication operation between the matrices and hence simple to implement, however, it was reported to have slow convergence \cite{slow_convergence}.
Gradient descent algorithms \cite{gradient1}, \cite{gradient2} are also employed to solve the NMF problem  in (\ref{eq:11}); wherein the ${\bw}$ and ${\bh}$ matrices are updated by taking a step in a direction opposite to the direction of gradient of the function $f_{_{\rm NMF}}\left(\bw, \bh\right)$. The update equation of a gradient descent algorithm for NMF problem is as follows:
\begin{equation}
\begin{array}{ll}
\bw^{k+1} = \bw^{k} - \alpha\dfrac{\partial f_{_{\rm NMF}}\left(\bw^{k},\bh^{k}\right)}{\partial \bw^{k}} = \bw^{k} - \alpha\left(\left(\bw^{k} \bh^{k}{\bh^{k}}^{T}\right)- \bv{\bh^{k}}^{T}\right)
\end{array}
\end{equation}
\begin{equation}
\begin{array}{ll}
\bh^{k+1} = \bh^{k} - \beta \dfrac{\partial f_{_{\rm NMF}}\left(\bw^{k+1},\bh^{k}\right)}{\partial \bh^{k}}= \bh^{k} - \beta\left(\left(\left({\bw^{k+1}}\right)^{T}\bw^{k+1}\bh^{k}\right)-{\left(\bw^{k+1}\right)}^{T}\bv\right)\\
\end{array}
\end{equation}
where $\alpha$ and $\beta$ are the step sizes. MU algorithm can also be viewed as the gradient descent algorithm \cite{mu_gradient} with iteration dependent or adaptive step size ${\alpha^{k}}$ and ${\beta^{k}}$ defined for updating matrices ${\bw}$ and ${\bh}$ respectively as: 
\begin{equation} \label{eq:12}
\begin{array}{ll}
\alpha^{k} =  \bw^{k} \oslash  {(\bw^{k} \bh^{k}{\bh^{k}}^{T})} \\
\beta^{k} = \bh^{k} \oslash  {\left({\left(\bw^{k+1}\right)}^{T}\bw^{k+1}\bh^{k}\right)}
\end{array}
\end{equation}
where ${\oslash}$ denotes element wise division.  Hence, the update equation becomes:
\begin{equation}\label{updatemu}
\begin{array}{ll}
\bw^{k+1} = \bw^{k}\circ \left(\left(\bv{\bh^{k}}^{T}\right) \oslash \left({\bw^{k}}\bh^{k}{\bh^{k}}^{T}\right)\right)\\
\bh^{k+1} = \bh^{k}\circ \left(\left({\bw^{k}}^{T}\bv\right) \oslash \left({\left(\bw^{k+1}\right)}^{T}\left(\bw^{k+1}\right){\bh^{k}}\right)\right)\\
\end{array}
\end{equation}
where ${\circ}$ denotes element wise multiplication. Note that the above update equation for $\bh$ and $\bw$ have only multiplication operation between the matrices. This kind of algorithm falls under multiplicative update algorithms. If in contrast, if the step sizes are chosen such that the update equation have only addition operations between the matrices, then it is called an additive update algorithm \cite{lecture_note}.
Another group of algorithms use the fact that steps 1 and 2 of Alternating minimization (refer to pseudocode of alternating minimization in Table 1) fall under \emph{Nonnegative Least Squares problem (NLS)} and solves each subproblem using active set methods (\cite{nlsa}, \cite{nlsb}, \cite{nlsc}), nevertheless, this approach can be computationally expensive \cite{gillis_thesis}. Instead of alternatingly minimizing over the two blocks - $\bw$ and $\bh$ matrices, one can form $2r$ blocks by further partitioning  the columns of $\bw$ matrix  and rows of $\bh$ matrix  as $\bw= $ $[\bwv_{1},\bwv_{2}\cdots \bwv_{r}]$, $\bh$ = $[\bhv_{1},\bhv_{2}\cdots \bhv_{r}]$ and solve the problem in (\ref{eq:11}) by updating the $\bwv_{j}^{th}$ and $\bhv_{j}^{th}$ block, while keeping the other blocks constant. Mathematically, this can be written as:
\begin{equation}\label{hals_a}
\begin{array}{ll}
\underset{\bwv_{j},\bhv_{j} \geq 0}{\rm minimize}\:{{\|\bv^{(j)}-\bwv_{j}\bhv_{j}^{T}\|}_{F}^{2}},  \textrm{for}\, \{j=1, 2 \cdots r\}
\end{array}
\end{equation}
where $\bv^{(j)} = \bv - \bw^{k}\bh^{k} -\bwv_{j}^{k}{\bhv_{j}^{k}}^{T}$. \emph{Hierarchical Alternating Least Square (HALS)} algorithm \cite{hals}, \cite{fasthals} follows this strategy and finds a closed form solution for the problem in (\ref{hals_a}) by alternatingly minimizing over $\bwv_{j}$ and $\bhv_{j}$. The advantage of this approach is that the computation of the solution for problem in (\ref{hals_a}) is computationally inexpensive compared to NLS. Fast-HALS is an efficient implementation of HALS algorithm wherein they do not explicitly compute $\bv^{(j)}$ and hence is computationally faster than HALS algorithm. Recently, Alternating Direction Method of Multipliers (ADMM) was also used to solve the NMF problem \cite{admm1}, \cite{admm2}. Vandaele et.al. \cite{new} proposed greedy randomized adaptive search procedure and simulated annealing to solve the NMF problem.  \\
\\ 
In this paper, we present two novel algorithms to solve the NMF problem-a sequential and a parallel algorithm. The major contributions of the paper are as follows:
\begin{enumerate}
\item{A Block MM based sequential algorithm-\textbf{I}terative \textbf{No}nnegative \textbf{M}atrix Factorization (\textbf{INOM}) is proposed to solve the NMF problem. The update equation of this algorithm looks like that of the update equation of gradient descent algorithm with iteration dependent or adaptive step size. Typically, many algorithms use line search to estimate the step size to be taken at every iteration but our algorithm is developed in such a way that the MM procedure gives us the step size to be taken at every iteration and hence we don't have to search for the optimal step size.}
\item{A parallel algorithm based on MM is presented, which parallely updates the $\bw$ and $\bh$ matrices. We address this algorithm as \textbf{Par}allel \textbf{I}terative \textbf{No}nnegative \textbf{M}atrix Factorization (\textbf{PARINOM}).}
\item{We discuss the convergence of the proposed algorithms and prove that they always converge to the stationary point.}
\item{Numerical simulations were conducted to compare the proposed algorithms with the existing algorithms and analyze the performance of the algorithms on the application of Blind Source Separation. }  
\end{enumerate}
The paper is organized as follows. An overview of MM and block MM can be found in section \ref{sec:2}. In section \ref{sec:3}, we first propose \textbf{INOM} to solve the NMF problem in (\ref{eq:11}). Next, we propose \textbf{PARINOM}  and in the same section, we also show that the proposed algorithms converge to stationary point and discuss the computational complexity of the algorithms. At the end of the section, we discuss an acceleration scheme to further accelerate the convergence of \textbf{PARINOM} algorithm. In section \ref{sec:6}, we compare the algorithms with the existing algorithms via computer simulations and evaluate the performance of our algorithm on the application of Blind Source Separation and conclude the paper in section \ref{sec:7}.\\
\\
Throughout the paper, {\bf{bold}} capital and {\bf{bold}} small letters are used to denote the matrix and vector, respectively. A scalar is denoted by a small letter. The ${i^{th}}$ entry of vector ${\bvv}$ is denoted by ${v_{i}}$ and ${(i,j)^{th}}$ entry of the matrix ${\bv}$ is denoted by ${v_{ij}}$. ${\circ}$ and ${\oslash}$ denotes element wise multiplication and division operation respectively. \\
\section{Majorization Minimization (MM): Vanilla MM and Block MM}\label{sec:2}
\subsection{Vanilla MM}
\begin{figure}[h]
\centering
\begin{tabular}{c}
\includegraphics[height=2.3in,width=3.8in]{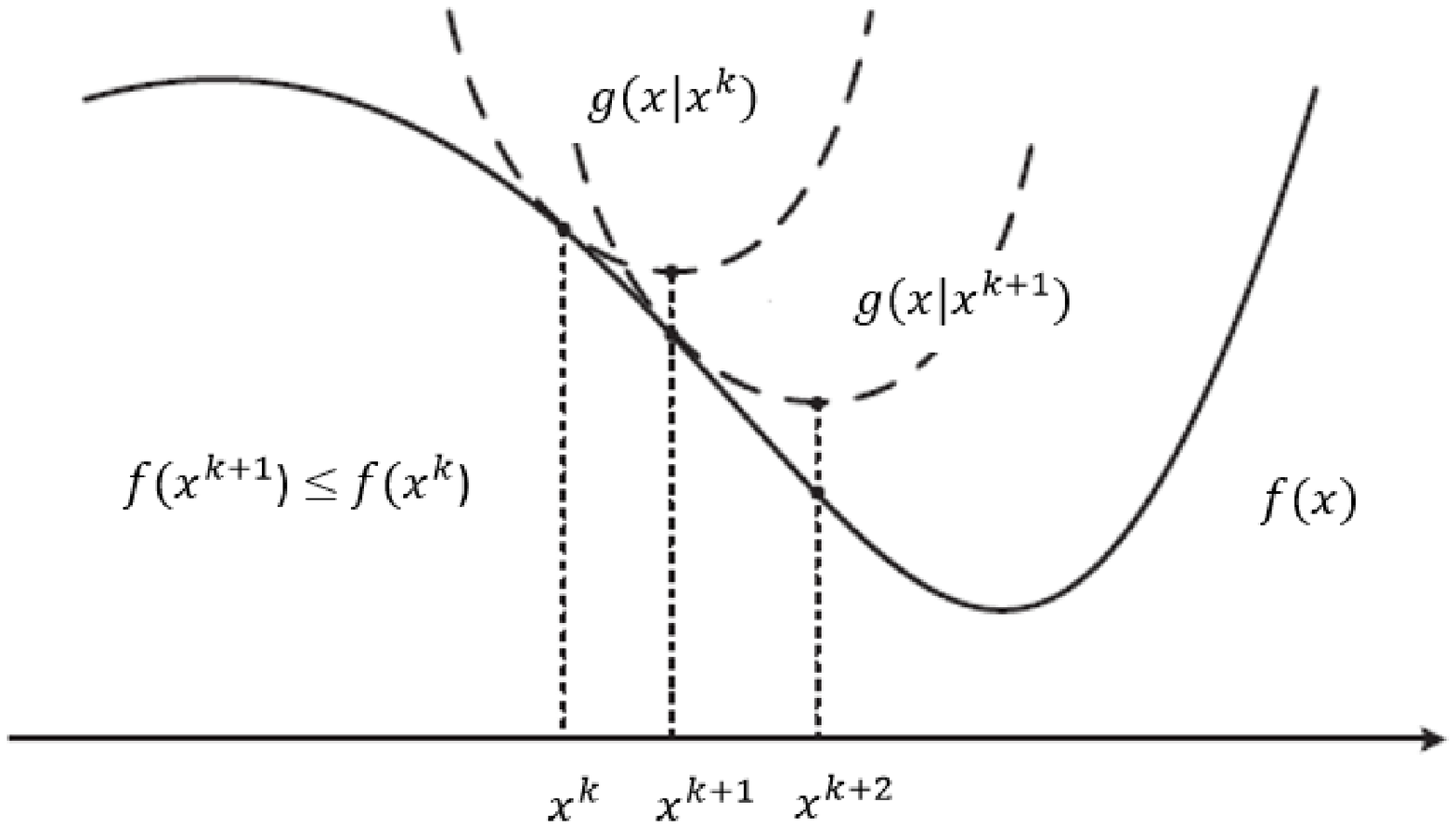}
\end{tabular}
\caption{MM procedure \cite{app_mm}}
\label{mm_procedure}
\end{figure}

Majorization Minimization (MM) is an iterative procedure which is mostly used to solve a non-convex, non-smooth or even a convex problem more efficiently. In MM, instead of minimizing the difficult optimization problem $f(\bxv)$ over the constraint set $\mathbf{\chi}$ directly, a ``surrogate'' function which  majorizes the problem (at a given value of $\bxv$ = $\bxv^{k} \in \chi$) is minimized at every iteration. The surrogate function $g(\bxv)$ is the global upper bound of the objective function $f(\bxv)$ i.e., it satisfies the following properties:
\begin{equation}  \label{eq:21}
g\left(\bxv^{k}|\bxv^{k}\right) = f\left(\bxv^{k}\right) 
\end{equation}
\begin{equation}\label{eq:22}
g\left(\bxv|\bxv^{k}\right) \geq f\left(\bxv\right)  \quad \textrm{for any}\, \bxv \in \chi
\end{equation}
where, ${\bxv^{k}}$ is the value taken by $\bxv$ at the $k^{th}$ iteration. Hence, the MM algorithm generates a sequence of points ${\{\bxv^{k}\}}$ according to the following rule: 
\begin{equation}  \label{eq:23}
\bxv^{k+1} \in \underset{\bxv \in \chi}{\rm arg\:min} \quad g\left(\bxv|\bxv^{k}\right)
\end{equation}
The MM procedure is depicted in Fig. \ref{mm_procedure}, wherein $g(\bxv|\bxv^{k})$ is the surrogate function which majorizes $f(\bxv)$ around $\bxv^{k}$ at the $k^{th}$ iteration. Pictorially, it can be seen that $f(\bxv^{k+2}) < f(\bxv^{k+1}) < f(\bxv^{k})$. 
By using (\ref{eq:21}), (\ref{eq:22}) and (\ref{eq:23}) it can be shown mathematically that the objective function is monotonically decreased at every iteration: \\
\begin{equation}\label{eq:24}
f(\bxv^{k+1}) \leq g\left(\bxv^{k+1}|\bxv^{k}\right) \leq g\left(\bxv^{k}|\bxv^{k}\right) = f(\bxv^{k})
\end{equation}
The first inequality and the last equality are by using (\ref{eq:21}) and (\ref{eq:22}), respectively. The second inequality is by (\ref{eq:23}). 
The convergence rate and computational complexity of the algorithm depends on how well one formulates the surrogate function. The convergence rate depends on how well the surrogate function follows the shape of the objective function $f(\bxv)$ and to have lower computational complexity, the surrogate function must be easy to minimize. Hence, the novelty of the algorithm based on MM lies in the design of the surrogate function.  Moreover, in the case of multivariate optimization problem, if the surrogate function is separable in the optimization variables - then the minimization problem could be solved parallely, which gives the computational advantage. To design surrogate function there are no set steps to follow. However, there are few papers which give guidelines for designing various surrogate functions \cite{tutorial}, \cite{app_mm}.
\\
\subsection{Block MM}
\begin{figure}[H]
\centering
\begin{tabular}{c}
\includegraphics[height=3in,width=6in]{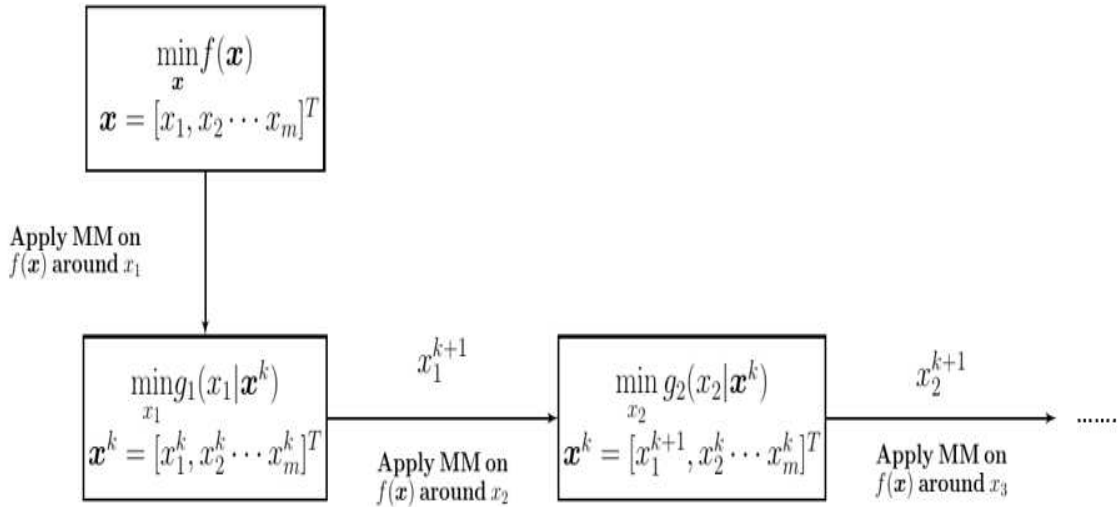}
\end{tabular}
\caption{Block MM procedure}
\label{block_mm_procedure}
\end{figure}

Suppose that the optimization variable ${\bxv}$ can be split into ${m}$ blocks as $\bxv =(\bxv_{1},\bxv_{2},\cdots,\bxv_{m})$, then one could apply Block MM to solve the optimization problem $f(\bxv)$. Block MM, an extension of vanilla MM, is a combination of block coordinate descent and vanilla MM wherein the optimization variable is split into several blocks and each block is updated using vanilla MM while keeping the other blocks fixed. Hence, the ${i}^{th}$ block is updated by minimizing the surrogate function $g_{i}\left(\bxv_{i}|\bxv^{k}\right)$ which majorizes $f(\bxv)$ on the ${i}^{th}$ block and has to satisfy the following properties:
\begin{equation}  \label{eq:24}
g_{i}{(\bxv_{i}\textsuperscript{k}|\bxv\textsuperscript{k})} = f(\bxv^{k}) 
\end{equation}
\begin{equation}\label{eq:25}
g_{i}(\bxv_{i}| \bxv^{k}) \geq f({\bxv_{1}}^{k},\cdots\bxv_{i},\cdots{\bxv_{m}}^{k}) 
\end{equation}
where $\bxv^{k}$ is the value taken by $\bxv$ at the $k^{th}$ iteration. The ${i}$th block at ${k+1}$ iteration is updated by solving the following problem:
\begin{equation}\label{s1}
{\bxv_{i}}^{k+1} \in \underset{\bxv_{i}}{\rm arg\:min} \quad g_{i}\left(\bxv_{i}|\bxv^{k}\right)
\end{equation}
Each block in Block MM is usually updated in a cyclic fashion. The Block MM procedure is as shown in Fig. \ref{block_mm_procedure}. 
The surrogate function for block MM must be chosen in such a way that the surrogate function is easy to minimize and must also approximate the objective function well to have faster convergence. Moreover, it is also reported that in some cases the surrogate function in case of block MM can approximate the objective function better than using a single block, leading to a faster convergence rate \cite{app_mm}. \\ 
\section{Algorithms for NMF problem}\label{sec:3}

In this section, we propose two iterative algorithms to solve the NMF problem in (\ref{eq:11}) based on block MM and vanilla MM. 
The first algorithm - \textbf{INOM} is based on block MM and sequentially updates both the matrices. The second algorithm - \textbf{PARINOM} is based on vanilla MM and parallely updates the $\bw$ and $\bh$ matrices. We then discuss the convergence of the two algorithms and show that they always converge to the stationary point. At the end of the section, we also show a way to accelerate the convergence of \textbf{PARINOM} algorithm, with a little additional computational cost. \\

\subsection{\textbf{INOM}:\textbf{I}terative \textbf{NO}nnegative \textbf{M}atrix Factorization}
This algorithm is based on Block Majorization Minimization principle wherein $\bw$ and $\bh$ are considered as blocks and each block is updated by vanilla MM scheme while keeping the other blocks fixed. \\
\\
Before we derive the update equation for $\bw$ and $\bh$ matrices, we first point out an important observation that is subsequently used to form the surrogate functions $g_{_{\bw}}\left(\bw|\bh^{k},\bw^{k}\right)$ and $g_{_{\bh}}\left(\bh|\bw^{k},\bh^{k}\right)$, which is used to majorize the function $f_{_{\rm NMF}}\left(\bw,\bh\right)$ on the $\bw^{th}$ and $\bh^{th}$ block, respectively. To discuss the same, we re-write $f_{_{\rm NMF}}\left(\bw,\bh\right)$  as: 
\begin{equation}\label{eq:35}
\begin{array}{ll}
f_{_{\rm NMF}}\left(\bw,\bh\right)= \|\bv - \bw\bh\|_{F}^{2}= \displaystyle\sum_{j=1}^{m}\|\bvv_{j} - \bw \bhv_{j}\|^{2}
\end{array}
\end{equation}
where $\bvv_{j}$ and $\bhv_{j}$ is the $j^{th}$ column of $\bv$ and $\bh$, respectively. From (\ref{eq:35}), it can be seen that given $\bw=\bw^{k}$, the $f_{_{\rm NMF}}\left(\bw^{k},\bh\right)$ is separable on each column of $\bh$. On taking the partial derivative of the above function with respect to ${\bhv_{j}}$, we get
\begin{equation}\label{eq:36}
\begin{array}{ll}
\nabla f_{_{\rm NMF}}(\bhv_ {j}|\bw^{k}) = -2(\bw^{k})^{T}\bvv_{j}+2(\bw^{k})^{T}\bw^{k}\bhv_{j} \\
\\
\nabla^{2} f_{_{\rm NMF}}(\bhv_{j}|\bw^{k}) = 2(\bw^{k})^{T}\bw^{k} 
\end{array}
\end{equation}
Note that the Hessian matrix $2(\bw^{k})^{T}\bw^{k}$ consists only of nonnegative entries and is independent of $\bhv_{j}$. We now show that $f_{_{\rm NMF}}\left(\bw,\bh^{k}\right)$ is separable on each row of $\bw$ and the Hessian of $f_{_{\rm NMF}}\left(\bw,\bh^{k}\right)$ with respect to $j^{th}$ row of $\bw$  is nonnegative and is equal to $2\bh^{k}(\bh^{k})^{T}$:
\begin{equation}\label{eq:1a}
\begin{array}{ll}
f_{_{\rm NMF}}\left(\bw,\bh^{k}\right)= \|\bv - \bw\bh^{k}\|_{F}^{2} = \displaystyle\sum_{j=1}^{n}\|\bvv'_{j} - \bwv_{j}\bh^{k}\|_{F}^{2}
\end{array}
\end{equation}
where $\bvv'_{j}$ and $\bwv_{j}$ represent the ${j^{th}}$ row of $\bv$ and $\bw$, respectively. 
\begin{equation}
\begin{array}{ll}
\nabla f_{_{\rm NMF}}(\bwv_ {j}|\bh^{k})=-2\bvv'_{j}{{\bh^{k}}}^{T} +2\bwv_{j}\bh^{k}{\bh^{k}}^{T}\\
\\
\nabla^{2} f_{_{\rm NMF}}(\bwv_{j}|\bh^{k}) = 2\bh^{k}{\bh^{k}}^{T}
\end{array}
\end{equation}
We smartly use the nonnegative Hessian matrix to design the upper bounds  $g_{_{\bw}}\left(\bw|\bh^{k},\bw^{k}\right)$ and $g_{_{\bh}}\left(\bh|\bw^{k},\bh^{k}\right)$, based on the following lemmas.

\begin{lemma} \label{lemma 1}
\begin{it}{ \bf {Lower Quadratic Bound Principle}} \end{it}\cite{tutorial}, \cite{quadratic_bound}: Given $\bxv=\bxv^{k}$, a twice differentiable function  $f(\bxv)$, a square matrix $\Lambda$ such that $\Lambda \succcurlyeq \nabla^{2} f(\bxv)$, then $f(\bxv)$ can be upper bounded as
\begin{equation}\label{eq:31} 
\begin{array}{ll}
f(\bxv) \leq f(\bxv^{k}) + \nabla f(\bxv^{k})^{T}(\bxv -\bxv^{k}) +\dfrac{1}{2}(\bxv -\bxv^{k})^{T}\Lambda(\bxv -\bxv^{k})
\end{array}
\end{equation}
The upper bound for $f(\bxv)$ is quadratic and differentiable in $\bxv$.
\end{lemma}
\begin{IEEEproof}
Suppose there exists a matrix $\Lambda$, such that $\Lambda \succcurlyeq \nabla^{2} f(\bxv)$, then we have the following equality by second order Taylor expansion:
\begin{equation*}
f(\bxv) \leq  f(\bxv^{k}) + \nabla f(\bxv^{k})^{T}(\bxv -\bxv^{k}) +\dfrac{1}{2}(\bxv -\bxv^{k})^{T}\Lambda(\bxv -\bxv^{k})
\end{equation*}
and equality is achieved at $\bxv =\bxv^{k}$. 
\end{IEEEproof}
\begin{lemma} \label{lemma 2}
Let ${\ba}$ of size $n \times n$ denote a nonnegative, square and symmetric matrix and let $p$ $\in$ $\mathbf{R^{+}}$. We define $\Lambda$ as a diagonal matrix with its diagonal elements equal to $p$, the maximum row sum of ${\ba}$ matrix. Then, $\Lambda \succcurlyeq {\ba}$.
\end{lemma}
\begin{IEEEproof}
The above lemma is based on {\bf{Perron-Frobenius Theorem}} \cite{pf_theorem1}, \cite{pf_theorem2}. According to the theorem, there is a positive real number $p$, called the Perron root, such that any other eigenvalue $\lambda$ of $\ba$ in absolute value is strictly smaller than $p$ i.e. $|\lambda| < p$.  To calculate $p$, an important corollary of the theorem is used:\\
\begin{equation}\label{eq:32}
\begin{array}{ll}
p \leq  \underset{i}{\rm max}\displaystyle\sum_{j=1}^{n} a_{ij}
\end{array}
\end{equation}
where $a_{ij}$ is the $(i,j)^{th}$ element of $\ba$ matrix. Hence, we get the following inequality:
\begin{equation}
\underset{i}{\rm max}\displaystyle\sum_{j=1}^{n} a_{ij} \geq p > \lambda
 \end{equation}
By constructing a diagonal matrix $\Lambda$ with its diagonal elements as maximum row sum of $\ba$ and by using Eigen Value Decomposition, the inequality $\Lambda \succcurlyeq {\ba}$ is attained. \\
\end{IEEEproof}
Using lemma \ref{lemma 1} and lemma \ref{lemma 2}, we construct the surrogate function $g_{_\bh}\left(\bh|\bh^{k},\bw^{k}\right)$ for the problem in  (\ref{eq:35}) to update the $\bh^{th}$ block with $\bw^{k}$ fixed.
\begin{equation}\label{surrogate1}
\begin{array}{ll}
g_{_\bh}\left(\bh|\bh^{k},\bw^{k}\right)= \displaystyle\sum_{j=1}^{m}g_{\bhv_{j}}\left({\bhv_{j}}|{\bhv_{j}}^{k},\bw^{k}\right)= \displaystyle\sum_{j=1}^{m} f(\bhv_{j}^{k}) + \nabla f(\bhv_{j}^{k})^{T}(\bhv_{j} -\bhv_{j}^{k}) +\dfrac{1}{2}(\bhv_{j} - \bhv_{j}^{k})^{T}\Lambda^{k}(\bhv_{j} - \bhv_{j}^{k})
\end{array}
\end{equation}
where $\Lambda^{k}$ matrix is a diagonal matrix with its diagonal elements as maximum row sum of ${2{(\bw^{k})}^{T}{\bw^{k}}}$ and $\nabla f(\bhv_{j}^{k})$ is equal to  $-2{(\bw^{k})}^{T}\bvv_{j}+2{(\bw^{k})}^{T}\bw^{k}{\bhv_{j}}^{k}$. The surrogate function is separable in each column of $\bh$. Hence at any iteration, given $\bw=\bw^{k}$ and $\bh=\bh^{k}$, the surrogate minimization problem is:
\begin{equation}\label{sp1}
\begin{array}{ll}
\underset{\bh \geq 0}{\rm minimize} \:g_{_\bh}\left(\bh|\bh^{k},\bw^{k}\right) 
\end{array}
\end{equation}
which has a closed form solution given by:
\begin{equation}\label{eq:37}
\begin{array}{ll}
\bhv_{j}^{k+1} = \bhv_{j}^{k} + \dfrac{1}{\mu^{k}}\left({2{{\left(\bw^{k}\right)}^{T}}{\bvv_{j}}-{2{{\left(\bw^{k}\right)}^{T}}{\bw^{k}}\bhv_{ j}^{k}}}\right) \\
\\
\bhv_{j}^{k+1} = \bhv_{j}^{k} + \dfrac{1}{\mu^{k}}\left({\left(2{{\left(\bw^{k}\right)}^{T}}\bv\right)_{j}} -{\left(2{\left(\bw^{k}\right)}^{T}\bw^{k}\bh^{k}\right)_{j}}\right)
\end{array}
\end{equation}
where $\mu^{k}$ = maximum row sum of ${2{(\bw^{k})}^{T}{\bw^{k}}}$. Since the update of $j^{th}$ column of $\bh$ does not depend on the update of $({j+1})^{th}$ column of $\bh$, we can re-write the above update equation as: 
\begin{equation}
\begin{array}{ll}\label{eq:38}
\bh^{k+1} = \bh^{k} + \dfrac{1}{\mu^{k}}\left(2{\left(\bw^{k}\right)}^{T}\bv - 2{\left(\bw^{k}\right)}^{T}{\bw^{k}}\bh^{k}\right)
\end{array}
\end{equation}

We now construct the surrogate function $g_{_\bw}\left(\bw|\bw^{k},\bh^{k+1}\right)$ using lemma \ref{lemma 1} and lemma \ref{lemma 2} for the problem in  (\ref{eq:1a}) to update the $\bw^{th}$ block with $\bh^{k+1}$ fixed. 
\begin{equation}\label{eq:1c}
\begin{array}{ll}
g_{_\bw}\left(\bw|\bw^{k},\bh^{k+1}\right)=\displaystyle\sum_{j=1}^{n}g_{\bwv_{j}}\left({\bwv_{j}}|{\bwv_{j}}^{k},\bh^{k}\right)=\displaystyle\sum_{j=1}^{n}f(\bwv_{j}^{k}) + (\bwv_{j}-\bwv_{j}^{k})\nabla f(\bwv_{j}^{k})^{T} +\dfrac{1}{2}(\bwv_{j}-\bwv_{j}^{k})\Lambda^{k}(\bwv_{j}-\bwv_{j}^{k})^{T}
\end{array}
\end{equation}
where ${\Lambda^{k}}$ matrix is constructed with its diagonal elements as maximum row sum of ${2(\bh^{k+1})({\bh^{k+1}})^{T}}$ and $\nabla f(\bwv_{j}^{k}) =-2\bvv'_{j}({{\bh}^{k+1}})^{T} +2\bwv_{j}(\bh^{k+1})({\bh^{k+1}})^{T} $. The surrogate function is separable in each row of $\bw$. Hence at any iteration, given $\bw = \bw^{k}$ and $\bh = \bh^{k+1}$, the surrogate minimization problem is: 
\begin{equation}\label{sp2}
\begin{array}{ll}
\underset{\bw \geq 0}{\rm minimize} \:g_{_\bw}\left(\bw|\bw^{k},\bh^{k+1}\right)
\end{array}
\end{equation}
which has a closed form solution given by:
\begin{equation}\label{eq:39}
\begin{array}{ll}
\bw^{k+1}  = \bw^{k} + \dfrac{1}{\nu^{k}}\left (2\bv{\bh^{k+1}}^{T}  - 2\bw^{k}\bh^{k+1}{\bh^{k+1}}^{T}\right)
\end{array}
\end{equation} 
where $\nu^{k}$ = maximum row sum of ${2\bh^{k+1}{\bh^{k+1}}^{T}}$. Note that once $\bh^{k+1}$ and $\bw^{k+1}$ are computed using (\ref{eq:38}) and (\ref{eq:39}), the negative elements must be projected back to $\mathbf{R^{+}}$, to satisfy the nonnegative constraint. Also, observe that the update equations looks similar to gradient descent algorithm with ``adaptive'' step sizes equal to $\mu^{k}$ and $\nu^{k}$. Typically, line search algorithm is implemented to  iteratively search for the optimal step sizes until the objective function decreases \cite{line_search}. Nevertheless, we don't have to iteratively search for the optimal step sizes; since the proposed algorithm is based on block MM-it is guaranteed that the step sizes - $\nu^{k}$ and $\mu^{k}$ will ensure the monotonic decrease of the respective cost functions. We now calculate the computational complexity of \textbf{INOM} algorithm. The complexity of computing $\mu^{k}$ and $\nu^{k}$ is $\mathcal{O}(r^{2}n)$ and $\mathcal{O}(r^{2}m)$, respectively. The complexity of computing $\bh^{k+1}$ is $\mathcal{O}(rnm)+ \mathcal{O}(r^{2}(n+m))$. The complexity of computing $\bw^{k+1}$ is $\mathcal{O}(rnm)+ \mathcal{O}(r^{2}(n+m))$. Since \textbf{INOM} is a sequential algorithm, the total complexity of the algorithm is  $\mathcal{O}(2rnm+2r^{2}(n+m))$. Pseudo code of \textbf{INOM} is shown in Table 2.
\begin{center}
\begin{tabular}{   l }
\hline
\hline
\textbf{Table 2: INOM} \\
\hline
\hline
{\bf{Input}}: Data samples ${\bv}$ with each column normalized; ${r}$, ${m}$ and ${n}$. \\
{\bf{Initialize}}: Set ${k}$ = 0. Initialize ${\bw^{0}}$ and ${\bh^{0}}$. Each column of ${\bw^{0}}$ is normalized. \\
{\bf{Repeat}}: \\
 Update ${\bh}$ \\
    1) $\mu^{k}$ = maximum row sum of ${2{\bw^{k}}^{T}{\bw^{k}}}$ \\
    2) $\bh^{k+1} = \bh^{k} + \dfrac{1}{\mu^{k}}\left(2{\bw^{k}}^{T}\bv - 2{\bw^{k}}^{T}{\bw^{k}}\bh^{k}\right)$\\
    3) $\bh^{k+1} = \textrm{max}(\bf{\bf{0}}, \bh^{k+1})$\\  
 Update ${\bw}$ \\ 
    4) $\nu^{k}$ = maximum row sum of ${2\bh^{k+1}{\bh^{k+1}}^{T}}$\\
    5) $\bw^{k+1}  = \bw ^{k} + \dfrac{1}{\nu^{k}}\left (2\bv{\bh^{k+1}}^{T}  - 2\bw^{k}\bh^{k+1}{\bh^{k+1}}^{T}\right)$\\   
    6) $\bw^{k+1} = \textrm{max}({\bf{0}}, \bw^{k+1})$\\   
    7) normalize each column of ${\bw}$\\
   $k \leftarrow k+1$,
  \bf{until} $\left|\dfrac{f_{_{\rm NMF}}\left(\bw^{k+1},\bh^{k+1}\right)-f_{_{\rm NMF}}\left(\bw^{k},\bh^{k}\right)}{f_{_{\rm NMF}}\left(\bw^{k},\bh^{k}\right)}\right| \leq 10^{-6}$ \\
\hline
\hline
\end{tabular}
\end{center}
\subsection{\textbf{PARINOM}:\textbf{Par}allel \textbf{I}terative \textbf{No}nnegative \textbf{M}atrix Factorization}
\textbf{PARINOM} solves the problem in (\ref{eq:11}) using Vanilla MM without alternatingly minimizing over $\bw$ and $\bh$ and hence at iteration $i+1$, $\bh^{i+1}$ matrix does not depend on $\bw^{i+1}$ matrix.  Therefore, $\bh^{i+1}$ and $\bw^{i+1}$ matrix can be parallely updated. On expanding $f_{_{\textrm{NMF}}}\left(\bw,\bh\right)$ we get:
\begin{equation} \label{eq:31.1}
\begin{array}{ll}
\|\bv -\bw\bh\|_F^2= {\rm Tr}(\bv^{T}\bv) + {\rm Tr}(\bh^{T}\bw^{T}\bw\bh) -2{\rm Tr}(\bv^{T}\bw\bh)
\end{array}
\end{equation}
Ignoring the constant terms in (\ref{eq:31.1}), we will now show that the second and third term are sigmoidal functions, which is defined as:\\
\emph{Sigmoidal function} \cite{sigmoidal_def}, \cite{sig_mm}: Let $c$ be a positive or a negative number and ${x_{1}, x_{2} \cdots x_{n}}$ be the nonnegative components of a $n$-dimensional vector ${\bxv}$. Let ${\alpha_{j}}$ be the ${j^{th}}$ component of {\boldmath{${\alpha}$}}, which is the fractional power (can be positive, negative or zero) of each component of ${\bxv}$. Then, $c\displaystyle\prod_{j=1}^{n} x_{j}^{\alpha_{j}}$ is called a sigmoidal function.\\
The second term can be re-written as:
\begin{equation} \label{eq:31.3}
\begin{array}{ll}
\displaystyle \sum_{k=1}^{m} \sum_{j =1}^{n} \left(\sum_{l=1}^{r}\sum_{m=1}^{r} w_{jl}h_{lk}w_{jm}h_{mk}\right)
\end{array}
\end{equation}
where $w_{ab}$ and $h_{ab}$ represent the $(a,b)^{th}$ element of $\bw$ and $\bh$ matrix. The terms in (\ref{eq:31.3}) are sigmoidal functions with positive coefficient. Now, we show that the third term in (\ref{eq:31.1}) can also be written as sigmoidal function. 
\begin{equation} \label{eq:31.5}
\begin{array}{ll}
- 2{\rm Tr}(\bv^{T}\bw\bh) =  -2 \displaystyle \sum_{j =1}^{n}\displaystyle \sum_{k =1}^{m} \left(v_{jk}\displaystyle \sum_{l =1}^{r} w_{jl}h_{lk}\right)
\end{array}
\end{equation}
This is a sigmoidal function with negative coefficient. Hence, the objective function $f_{_{\textrm{NMF}}}\left(\bw,\bh\right)$ becomes:
\begin{equation}\label{new_obj}
\begin{array}{ll}
\|\bv -\bw\bh\|_F^2 = \displaystyle \sum_{k=1}^{m} \sum_{j =1}^{n} \left(\sum_{l=1}^{r}\sum_{m=1}^{r} w_{jl}h_{lk}w_{jm}h_{mk}\right) -2 \displaystyle \sum_{j =1}^{n}\displaystyle \sum_{k =1}^{m} \left(v_{jk}\displaystyle \sum_{l =1}^{r} w_{jl}h_{lk}\right)
\end{array}
\end{equation}
Now, we propose the following lemma which is used to majorize  the above objective function. 

\begin{lemma} \label{lemma 3} When c in $c\displaystyle\prod_{j = 1}^{n} x_{j}^{\alpha_{j}}$ is positive, the following inequality holds:
\begin{equation}\label{eq: 33}
\begin{array}{ll}
\textrm{c}\displaystyle \prod_{j=1}^{n} x_{j}^{\alpha_{j}} \leq \textrm{c}\displaystyle\prod_{j=1}^{n} \left({{x_{j}}^{i}}\right)^{\alpha_{j}}\displaystyle\sum_{j=1}^{n} \dfrac{\alpha_{j}}{\|{\alpha}\|_{1}}\left(\dfrac{x_{j}}{x_{j}^{i}}\right)^{\|{\alpha}\|_{1}}
\end{array}
\end{equation}
The above inequality is equal when ${\bxv}$ is equal to ${\bxv^{i}}$. When c is negative, the inequality in (\ref{eq: 33}) changes to:
\begin{equation}\label{eq: 34}
\begin{array}{ll}
\textrm{c}\displaystyle \prod_{j = 1}^{n} x_{j}^{\alpha_{j}} \leq \textrm{c}\left(\displaystyle\prod_{j = 1}^{n} \left({{x_{j}}^{i}}\right)^{\alpha_{j}}\right) \alpha_{j} \textrm{ln}(x_{j})
\end{array}
\end{equation}
\end{lemma}
\begin{IEEEproof}
See [Section 3, \cite{sig_mm}].\\
\end{IEEEproof}
Using the above lemma we get the following surrogate function $g\left(w_{jl},h_{lk}|w_{jl}^{i},h_{lk}^{i}\right)$ for $f_{_{\textrm{NMF}}}\left(\bw,\bh\right)$.
\begin{equation} \label{eq:31.8}
\begin{array}{ll}
g\left(w_{jl},h_{lk}|w_{jl}^{i},h_{lk}^{i}\right) = \displaystyle\sum_{k=1}^{m}  \displaystyle\sum_{j =1}^{n} \left( \displaystyle\sum_{l=1}^{r} \displaystyle\sum_{m=1}^{r} w_{jl}^{i}h_{lk}^{i}w_{jm}^{i}h_{mk}^{i} \left(\dfrac {1}{4} \left(\dfrac{w_{jl}}{w_{jl}^{i}}\right)^4 + \dfrac {1}{4} \left(\dfrac{h_{lk}}{h_{lk}^{i}}\right)^4  +\dfrac {1}{4} \left(\dfrac{w_{jm}}{w_{jm}^{i}}\right)^4 + \dfrac {1}{4} \left(\dfrac{h_{mk}}{h_{mk}^{i}}\right)^4 \right)\right)   \\
-2\displaystyle\sum_{j =1}^{n} \displaystyle\sum_{k =1}^{m}\left(v_{jk} \displaystyle\sum_{l =1}^{r}\left(w_{jl}^{i} h_{lk}^{i}\right)\left(\textrm{ln}w_{jl} + \textrm{ln} h_{lk}\right)\right)
\end{array}
\end{equation}
where $w_{ab}^{i}$ and $h_{ab}^{i}$ represent the $(a,b)^{th}$ element of $\bw$ and $\bh$ matrix at the $i^{th}$ iteration. Note that the surrogate function 
$g\left(w_{jl},h_{lk}|w_{jl}^{i},h_{lk}^{i}\right)$ is separable in the optimization variables:  $w_{jl}$ and $h_{lk}$.  Hence at any iteration, given $w_{jl} = w_{jl}^{i}$ and $h_{lk} = h_{lk}^{i}$, the surrogate minimization problem is:
\begin{equation}
\begin{array}{ll}
\underset{w_{jl} >0,\,h_{lk}>0}{\rm minimize} \: g\left(w_{jl}, h_{lk}|w_{jl}^{i},h_{lk}^{i}\right)
\end{array}
\end{equation}
where $g\left(w_{jl},h_{lk}|w_{jl}^{i},h_{lk}^{i}\right)$ is given by (\ref{eq:31.8}) which has a closed form solution is given by:
\begin{equation} \label{eq:31.9}
\begin{array}{ll}
w_{jl}^{i+1} =\sqrt[4]{\dfrac{{p_{(jl)}^{i}}(w_{jl}^{i})^4}{z_{1(jl)}^{i}}} \\
\\
h_{lk}^{i+1} =\sqrt[4]{\dfrac{{q_{(lk)}^{i}}(h_{lk}^{i})^4}{z_{2(lk)}^{i}}}
\end{array}
\end{equation}
where
\begin{equation} \label{eq:32.1}
\begin{array}{ll}
z_{1(jl)}^{i}= \displaystyle\sum_{k=1}^{m}\displaystyle\sum_{M=1}^{r} h_{lk}^{i}w_{jM}^{i}h_{Mk}^{i} \quad p_{(jl)}^{i}=  \displaystyle\sum_{k=1}^{m} v_{jk}h_{lk}^{i}\\
\\
z_{2(lk)}^{k}= \displaystyle\sum_{j=1}^{n} \displaystyle\sum_{M=1}^{r}  w_{jl}^{i}w_{jM}^{k}h_{Mk}^{i} \quad q_{(lk)}^{i} =\displaystyle\sum_{j=1}^{n} v_{jk}w_{jl}^{i}
\end{array}
\end{equation}
Taking the entire matrix into consideration, (\ref{eq:31.9}) can be re-written as: 
\begin{equation} \label{eq:32.2}
\begin{array}{ll}
\bw^{i+1} =\sqrt[4]{\left({ \left(\bv{\bh^{i}}^{T}\right)\circ {\bw^{i}}^4}\right)\oslash\left(\bw^{i}\bh^{i}{\bh^{i}}^{T}\right)} \\
\bh^{i+1} =\sqrt[4]{\left({ \left({\bw^{i}}^{T}\bv\right)\circ {\bh^{i}}^4}\right)\oslash\left({\bw^{i}}^{T}\bw^{i}\bh^{i}\right)}
\end{array}
\end{equation}
In (\ref{eq:32.2}), the fourth power of the matrix $\bw$ and $\bh$ are done element wise. When compared to MU and \textbf{INOM} algorithm, in \textbf{PARINOM}, update of $\bh^{i+1}$ does not depend on $\bw^{i+1}$ at ${i+1}$ iteration. Hence, $\bw^{i+1}$ and $\bh^{i+1}$ matrices can be updated parallely. The complexity in computing  $\bw^{i+1}$ is $\mathcal{O}(nmr + r(m+n))$. The complexity to update ${\bh^{i+1}}$ is also $\mathcal{O}(nmr +r(m+n))$. Since both the matrices can be updated parallely, the complexity of \textbf{PARINOM} is $\mathcal{O}(nmr + r(m+n))$. The Pseudo code of \textbf{PARINOM} is shown in Table 3:
\begin{center}
\begin{tabular}{   l }
\hline
\hline
\bf{Table 3: PARINOM} \\
\hline
\hline
{\bf{Input}}: Data samples ${\bv}$, with each column normalized, $r$, $m$ and $n$. \\
{\bf{Initialize}}: Set \emph{i} = 0. Initialize ${\bw^{0}}$ and ${\bh^{0}}$. Each column of ${\bw^{0}}$ is normalized. \\
{\bf{Repeat}}: \\
    Update ${\bh}$  and ${\bw}$ {\bf{parallely}} \\
    1)$\bw^{i+1} =\sqrt[4]{\left({ \left(\bv{\bh^{i}}^{T}\right)\circ {\bw^{i}}^4}\right)\oslash\left(\bw^{i}\bh^{i}{\bh^{i}}^{T}\right)}$  \\
    2)$\bh^{i+1} =\sqrt[4]{\left({ \left({\bw^{i}}^{T}\bv\right)\circ {\bh^{i}}^4}\right)\oslash\left({\bw^{i}}^{T}\bw^{i}\bh^{i}\right)}$  \\
    normalize each column of ${\bw^{i+1}}$\\
    $i \leftarrow i+1$, \bf{until} $\left|\dfrac{f_{_{\textrm{NMF}}}\left(\bw^{i+1},\bh^{i+1}\right)-f_{_{\textrm{NMF}}}\left(\bw^{i},\bh^{i}\right)}{f_{_{\textrm{NMF}}}\left(\bw^{i},\bh^{i}\right)}\right| \leq 10^{-6}$ \\
\\
\hline
\hline
\end{tabular}
\end{center}

\subsection{Convergence Analysis}\label{sec:4}

We first discuss the convergence of \textbf{INOM}, which is based on Block MM. To discuss the same, we first prove that the surrogate functions $g_{_\bh}\left(\bh|\bh^{k},\bw^{k}\right)$ in (\ref{surrogate1}) and $g_{_\bw}\left(\bw|\bw^{k},\bh^{k+1}\right)$ in (\ref{eq:1c}) are quasi-convex and the problems in (\ref{sp1}) and (\ref{sp2}) have a unique minimum
\begin{IEEEproof}
From (\ref{surrogate1}) and (\ref{eq:1c}), it can be seen that $g_{_\bh}\left(\bh|\bh^{k},\bw^{k}\right)$ and $g_{_\bw}\left(\bw|\bw^{k},\bh^{k+1}\right)$ are separable in the columns and rows of $\bh$ and $\bw$, respectively i.e. $g_{_\bh}\left(\bh|\bh^{k},\bw^{k}\right) = \displaystyle\sum_{j=1}^{m}g_{\bhv_{j}}\left({\bhv_{j}}|{\bhv_{j}}^{k},\bw^{k}\right)$ and $g_{_\bw}\left(\bw|\bw^{k},\bh^{k+1}\right) = \displaystyle\sum_{j=1}^{n}g_{\bwv_{j}}\left({\bwv_{j}}|{\bwv_{j}}^{k},\bh^{k}\right)$. The Hessian of  $g_{\bhv_{j}}\left({\bhv_{j}}|{\bhv_{j}}^{k},\bw^{k}\right)$ and $g_{\bwv_{j}}\left({\bwv_{j}}|{\bwv_{j}}^{k},\bh^{k}\right)$  for every $j$ is $\Lambda$ matrix - which is a diagonal matrix made of nonnegative elements and hence is positive semi-definite. This implies that $g_{\bhv_{j}}\left({\bhv_{j}}|{\bhv_{j}}^{k},\bw^{k}\right)$ and $g_{\bwv_{j}}\left({\bwv_{j}}|{\bwv_{j}}^{k},\bh^{k}\right)$ are convex functions. Since the sum of convex functions is convex, $g_{_\bh}\left(\bh|\bh^{k},\bw^{k}\right)$ and $g_{_\bw}\left(\bw|\bw^{k},\bh^{k+1}\right)$ are also convex functions. Since every convex function has convex sublevel sets \cite{convex_book}, $g_{_\bh}\left(\bh|\bh^{k},\bw^{k}\right)$ and $g_{_\bw}\left(\bw|\bw^{k},\bh^{k+1}\right)$ are also quasi-convex functions. Also, the problems in (\ref{sp1}) and (\ref{sp2}) has a unique minimum; since we are minimizing a convex function. 
\end{IEEEproof}
Razaviyayn et. al. in Theorem 2.(a) of \cite{convergence} showed that the sequence of points generated by Block MM converge to the stationary point, provided the surrogate function is quasi-convex and the minimizer of the surrogate minimization problem is unique. Hence, \textbf{INOM} algorithm converges to the stationary point of the problem in (\ref{eq:11}) which is a direct application of Theorem 2. (a) in \cite{convergence}. \\
\\
We now discuss the convergence of \textbf{PARINOM}; which is based on Vanilla MM. Note that $f_{_{\rm NMF}}\left(\bw,\bh\right)$ in (\ref{eq:11}) is bounded below by zero and the constraint set is closed and convex. Also from (\ref{eq:24}), the sequence of points $\{\bw^{k},\bh^{k}\}$  monotonically decrease the NMF problem. Hence, the sequence $f_{_{\rm{NMF}}}(\bw^{k},\bh^{k})$ generated by \textbf{PARINOM} will at the least converge to the finite value. \\
\\
We now show that the sequence $\{\bw^{k},\bh^{k}\}$ will converge to the stationary point. To prove the same, we first group the variables $\bw$, $\bh$ into a single block $\bX$.  
From (\ref{eq:24}), we have:
\begin{equation}\label{conv1}
\begin{array}{ll}
f_{_{\rm NMF}}\left(\bX^{0}\right) \geq f_{_{\rm NMF}}\left(\bX^{1}\right) \geq f_{_{\rm NMF}}\left(\bX^{2}\right) \cdots
\end{array}
\end{equation}
Assume that there is a subsequence ${\bX^{r_{j}}}$ converging to a limit point $\bz$. Then, from (\ref{eq:21}), (\ref{eq:22}) and from (\ref{conv1}) we obtain:
\begin{equation}
\begin{array}{ll}
g\left(\bX^{r_{j+1}}|\bX^{r_{j+1}}\right)= f_{_{\rm NMF}}\left(\bX^{r_{j+1}}\right) \leq f_{_{\rm NMF}}\left(\bX^{r_{j}+1}\right) \leq g\left(\bX^{r_{j}+1}|\bX^{r_{j}}\right)\leq g\left(\bX|\bX^{r_{j}}\right)
\end{array}
\end{equation}
Letting \emph{j} $\rightarrow$ $\infty$, we get\\
\begin{equation}
\begin{array}{ll}
g\left(\bz|\bz\right)\leq g\left(\bX|\bz\right)
\end{array}
\end{equation}
which implies $g'(\bz|\bz) \geq 0$. Since the first order behavior of surrogate function is same as function $f\left(\cdot\right)$, (\cite{convergence}), $g'(\bz|\bz) \geq 0$ implies $f_{_{\rm NMF}}'(\bz) \geq 0$. Hence, $\bz$ is the stationary point of $f_{_{\rm NMF}}\left(\cdot\right)$ and therefore the proposed algorithm converges to the stationary point.
\subsection{Squarem Acceleration Scheme}\label{sec:5}
We now describe a way to further accelerate the convergence of \textbf{PARINOM} algorithm based on  Squared Iterative method (SQUAREM) \cite{acc_mm} acceleration scheme. Originally, this scheme was proposed for fixed - point Expectation Maximization algorithm. However, since MM is a generalization of EM, this scheme could be used for accelerating MM based algorithms as well. SQUAREM is based on Cauchy - Brazilai - Brownein method, which is a combination of Cauchy and BB (Brazilai and Browein) method to accelerate the convergence rate. This acceleration scheme is proven to be monotonic. \\
\\
Let ${{f_{1}(\bw_{k},\bh_{k})}}$ and ${{f_{2}(\bw_{k},\bh_{k})}}$ denote fixed - point functions which update $\bw$ and $\bh$ respectively i.e 
\begin{equation}\label{eq:41}
\begin{array}{ll}
\bw^{k+1}= f_{1}\left(\bw^{k},\bh^{k}\right) =  \sqrt[4]{({ \bv(\bh^{k})^{T}\circ (\bw^{k})^4})\oslash(\bv^{k}((\bh^{k})^{T}))} \\ 
\bh^{k+1} = f_{2}\left(\bw^{k},\bh^{k}\right)=  \sqrt[4]{({(\bw^{k})^{T}\bv\circ (\bh^{k})^4})\oslash((\bw^{k})^{T}(\bv^{k}))}
\end{array}
\end{equation}
The pseudo code of the acceleration scheme is shown in Table 4.
\begin{center}
\begin{tabularx}{\textwidth}{@{}lX@{}}
\hline
\hline
\bf{Table 4: Acc-PARINOM} \\
\hline
\hline
\endhead
{\bf{Input}}: Data samples ${\bv}$, with each column normalized, $r$, $m$ and $n$. \\
{\bf{Initialize}}: Set ${k}$ = 0. Initialize ${\bw^{0}}$ and ${\bh^{0}}$. Each column of ${\bw^{0}}$ is normalized. \\
{\bf{Repeat}}: \\
   
    1. Parallely update: $\bw^{1}= f_{1}(\bw^{k},\bh^{k}), \bh^{1} = f_{2}(\bw^{k},\bh^{k})$ \\
    2. Normalize each column of ${\bw^{1}}$ \\
    3. Parallely update: $\bw^{2}= f_{1}(\bw^{1},\bh^{1}), \bh^{2} = f_{2}(\bw^{1},\bh^{1})$  \\
    4. Normalize each column of ${\bw^{2}}$ \\
    5. Compute: ${r_{h} = \bh^{1} - \bh^{k}}$, ${v_{h} =  \bh^{2} - \bh^{1} - r_{h}}$ and ${\alpha_{h} = -\dfrac{\|r_{h}\|}{\|v_{h}\|}}$\\
    6. Parallely compute $r_{w},v_{w},\alpha_{w}$ \\
    7. $\bh = \textrm{max}\left(0, \bh^{k} - 2\alpha_{h}r_{h} +\alpha_{h}^{2}v_{h}\right)$\\
    8. $\bw = \textrm{max}\left(0, \bw^{k} - 2\alpha_{w}r_{w} +\alpha_{w}^{2}v_{w}\right)$\\
    9. Normalize each column of ${\bw}$ \\
    10. {\bf{while}}$\, {{\|\bv-\bw\bh\|}_{F} > {\|\bv - \bw^{k}\bh^{k}\|}_{F} }$ {\bf{do}}\\
    11.  \quad $\alpha_{h} \leftarrow \dfrac{\alpha_{h} -1}{2}$ , $\alpha_{w} \leftarrow \dfrac{\alpha_{w} -1}{2}$ \\
    12. \quad $\bh = \textrm{max}\left(0, \bh^{k} - 2\alpha_{h}r_{h} +\alpha_{h}^{2}v_{h}\right)$\\
    13. \quad $\bw = \textrm{max}\left(0, \bw^{k} - 2\alpha_{w}r_{w} +\alpha_{w}^{2}v_{w}\right)$\\
    14. \bf{end while}\\
    15. ${\bw^{k+1}} =\bw$, ${\bh^{k+1}} = \bh$ \\ 
    16. normalize each column of ${\bw^{k+1}}$\\
    17. $k \leftarrow k+1$\\
{\bf{until convergence}}\\
\hline
\hline
\end{tabularx}
\end{center}
SQUAREM acceleration scheme will sometimes violate the non-negative constraint. Hence, the negative values of the matrix must be made zero. To retain the descent property of the proposed algorithm; the value of $\alpha_{h}$ and $\alpha_{w}$ is found by backtracking which halves the distance between alpha and -1 until the descent property is satisfied. Note that when $\alpha_h$ is equal to -1, $\textrm{max}\left(0, \bh^{k} - 2\alpha_{h}r_{h} +\alpha_{h}^{2}v_{h}\right)$ becomes equal to ${\bh^{2}}$. Due to the monotonic behavior of MM, ${\bh^{2}}$ will be less than ${\bh^{k}}$. Similarly when $\alpha_w$ is equal to -1, $\textrm{max}\left(0, \bw^{k} - 2\alpha_{w}r_{w} +\alpha_{w}^{2}v_{w}\right)$ is equal to ${\bw^{2}}$ - which is lesser than ${\bw^{k}}$. Hence the descent property is guaranteed to hold as alpha is pushed towards -1. \\

\section{Simulations}\label{sec:6}
In this section, we present numerical simulations to compare the proposed methods with the state of the art algorithms - Fast-HALS and MU. To have a fair comparison we accelerate the MU algorithm using the SQUAREM acceleration scheme. All the simulations were carried out on a PC with 2.40GHz Intel Xeon Processor with 64 GB RAM. 

\begin{figure}[H]
\centering
\begin{subfigure}{0.49\textwidth}
\centering
\captionsetup{justification=centering}
\includegraphics[height=2.5in,width=3.3in]{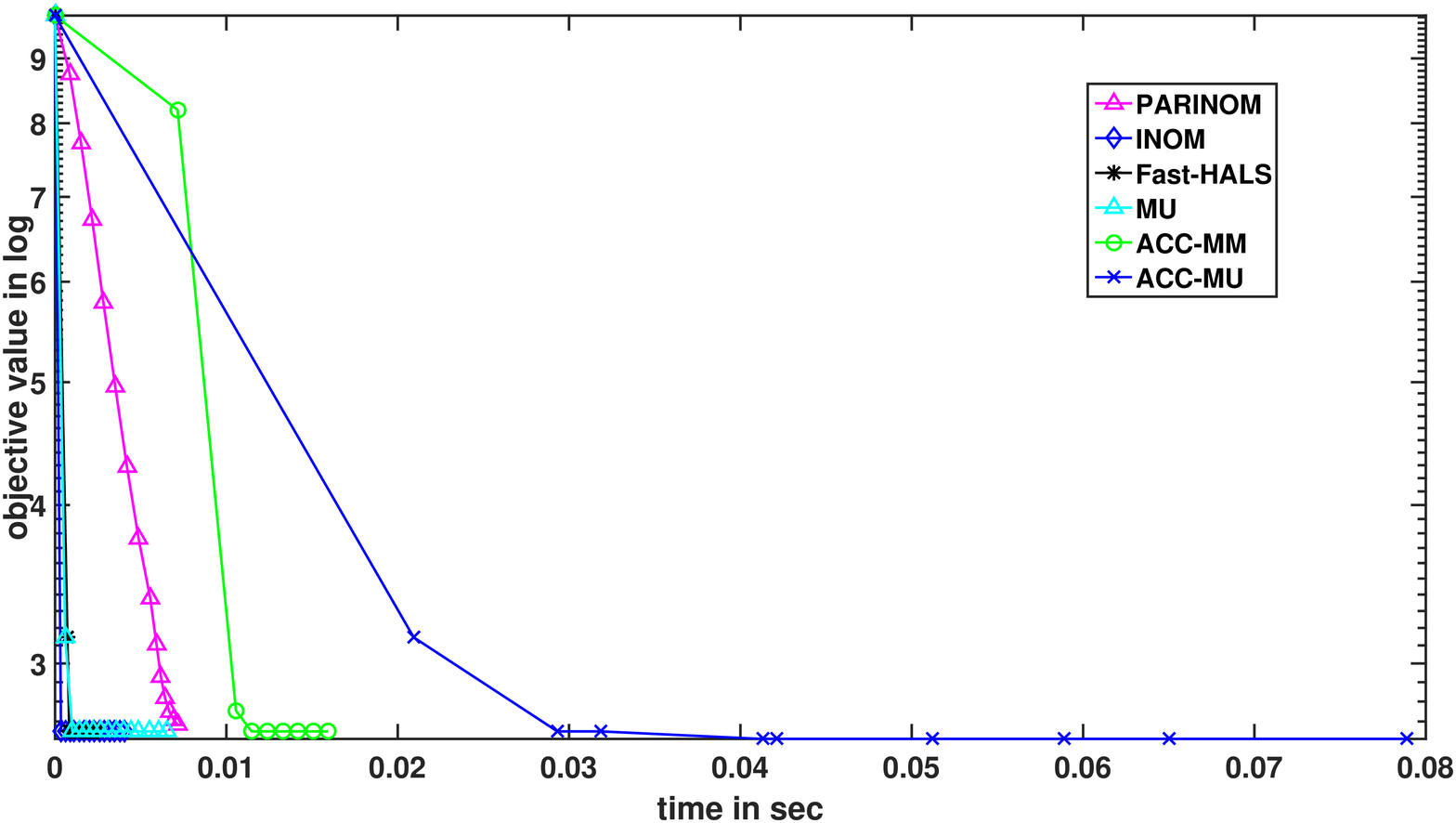}
\caption{Objective value vs run time of the proposed and existing algorithms- Fast-HALS, MU and Accelerated MU}
\end{subfigure}
\begin{subfigure}{0.49\textwidth}
\centering
\captionsetup{justification=centering}
\includegraphics[height=2.5in,width=3.3in]{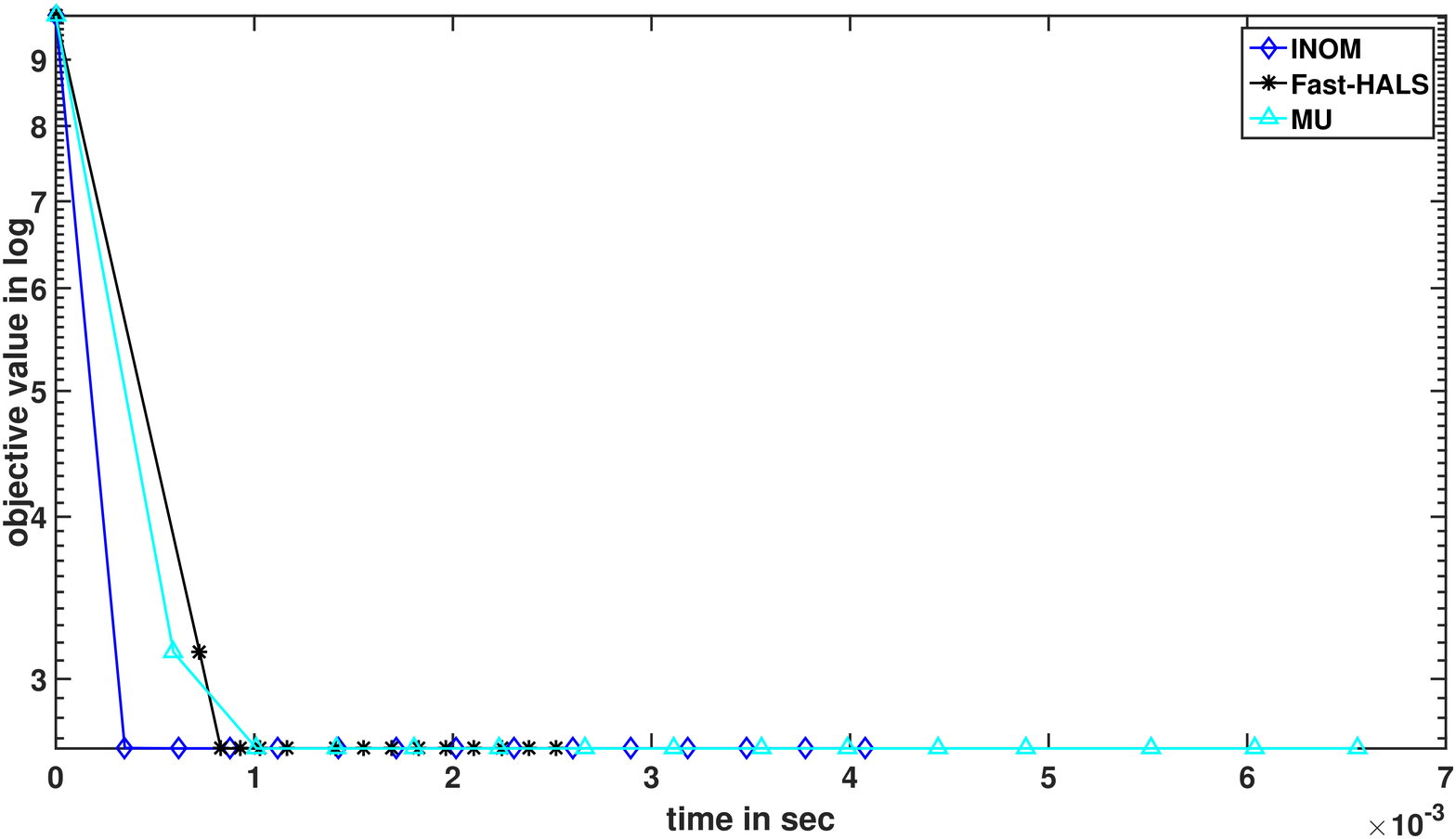}
\caption{Objective value vs run time of MU,\\ Fast-HALS  and \textbf{INOM}}
\end{subfigure}
\caption{Comparison of run time and convergence of proposed algorithm with existing algorithms when $n=100$, $m=200$ and $r=1$.}
\label{fig1}
\end{figure}

1. In the first simulation, we fix $n=100$, $m=200$ and $r=1$ and compare the convergence rate and per iteration cost of the proposed algorithms vs Fast-HALS, MU and Accelerated MU. The elements of $\bv$ matrix was randomly generated from a uniform distribution from $[100, 200]$. Initial values of $\bw$ and $\bh$ was also randomly generated from a uniform distribution from $[0, 1]$. The columns of $\bv$ and $\bw$ were normalized. The initial objective value $f_{_{\textrm{NMF}}}(\bw^{0},\bh^{0})$ for all the algorithms were kept same. Fig. \ref{fig1} (a) shows the run time vs the objective value in log scale of all the algorithms. Fig. \ref{fig1} (b) compares the run time vs the objective value in log scale of MU, Fast-HALS and \textbf{INOM}.

\begin{figure}[H]
\centering
\begin{subfigure}{0.49\textwidth}
\centering
\captionsetup{justification=centering}
\includegraphics[height=2.5in,width=3.3in]{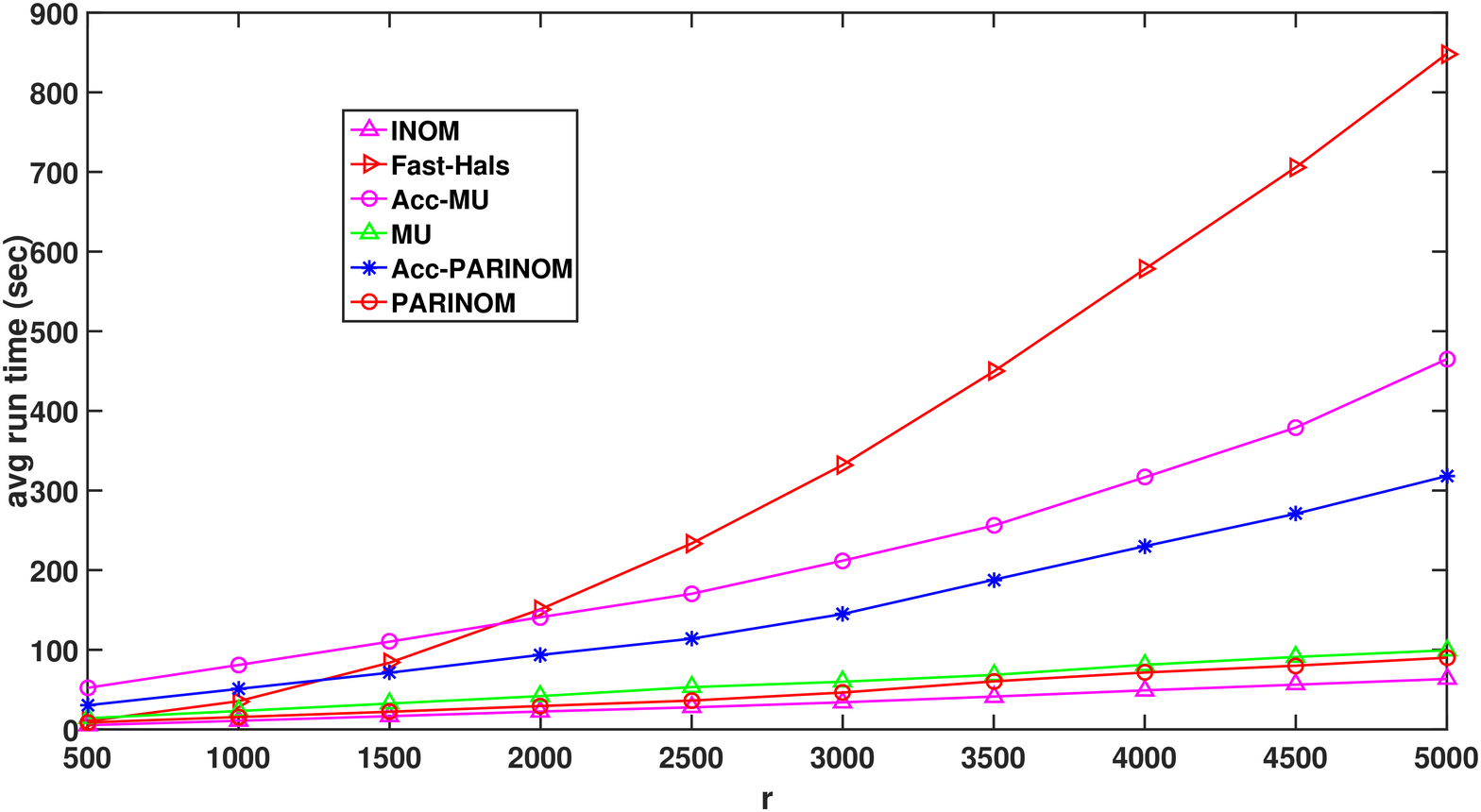}
\caption{Comparison of run time of proposed algorithms with Fast-Hals, MU and Acc-MU.}
\end{subfigure}
\begin{subfigure}{0.49\textwidth}
\centering
\captionsetup{justification=centering}
\includegraphics[height=2.5in,width=3.3in]{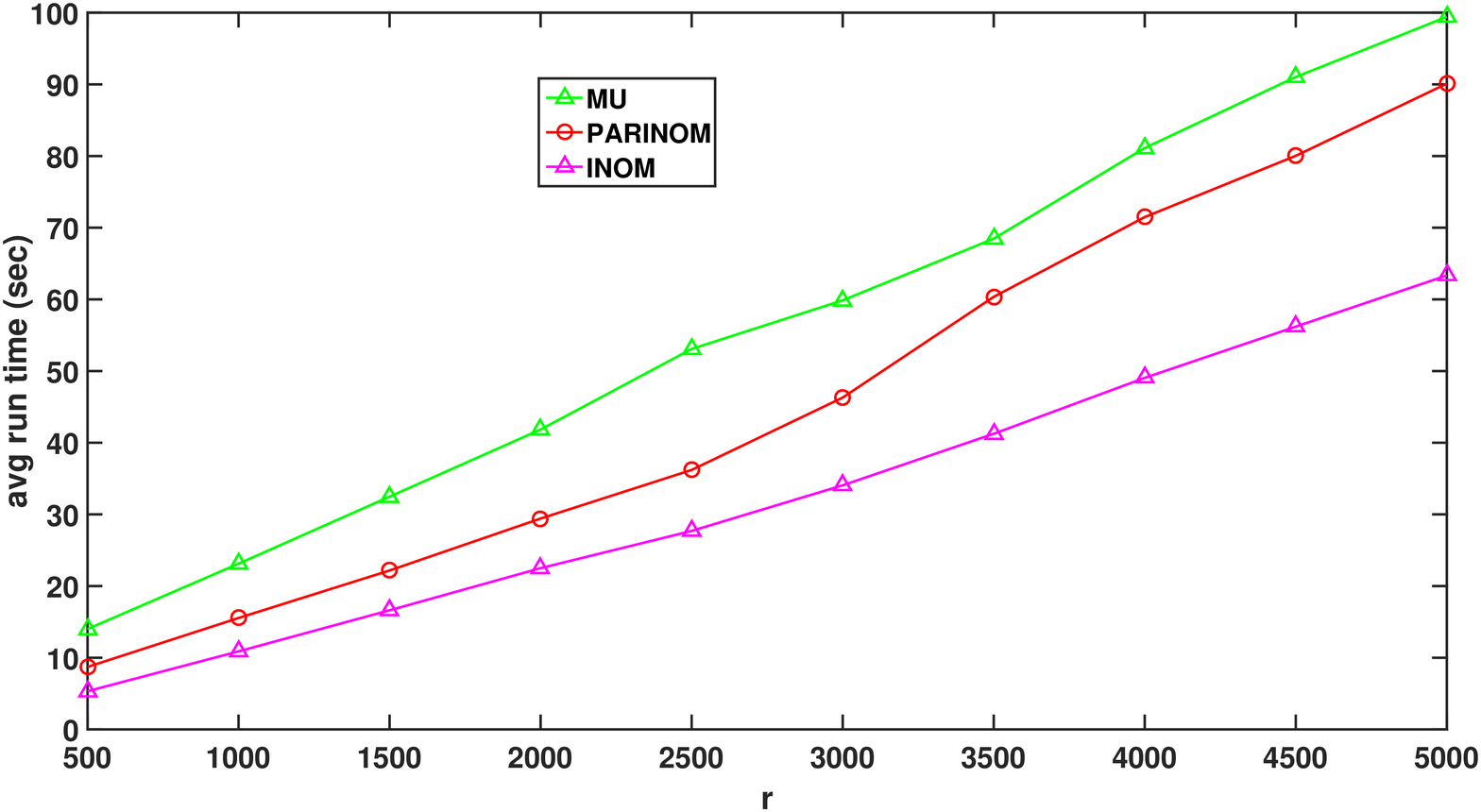}
\caption{Comparison of run time of MU, INOM and PARINOM}
\end{subfigure}
\caption{Comparison of proposed algorithm with existing algorithms by varying $r$, $m=50000$ and $n=10000$ and when $\bv$ is a dense matrix.}
\label{large_matrix}
\end{figure}

\begin{figure}[H]
\centering
\begin{subfigure}{0.49\textwidth}
\centering
\captionsetup{justification=centering}
\includegraphics[height=2.5in,width=3.3in]{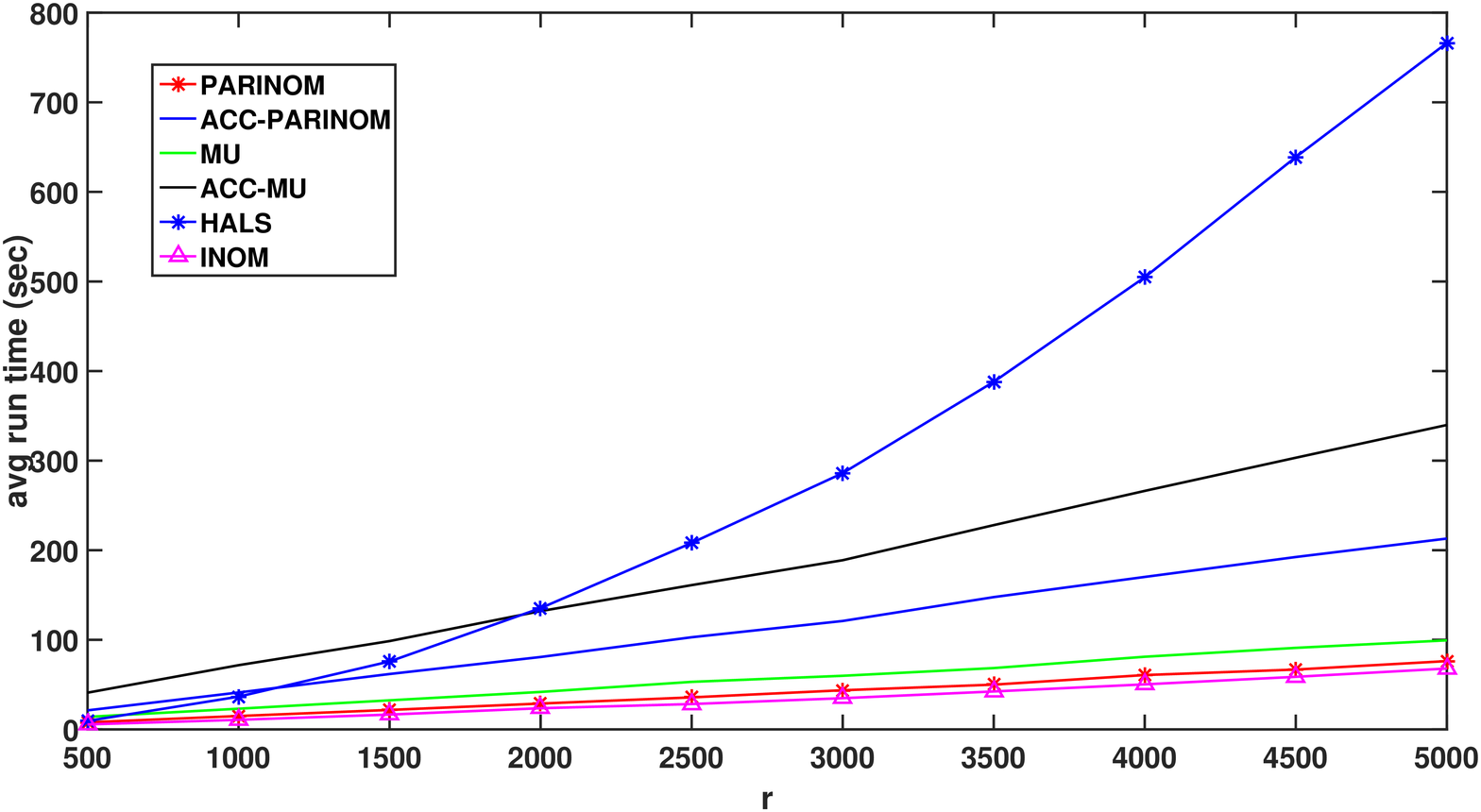}
\caption{Comparison of run time of proposed algorithms with Fast-Hals, MU and Acc-MU.}
\end{subfigure}
\begin{subfigure}{0.49\textwidth}
\centering
\captionsetup{justification=centering}
\includegraphics[height=2.5in,width=3.3in]{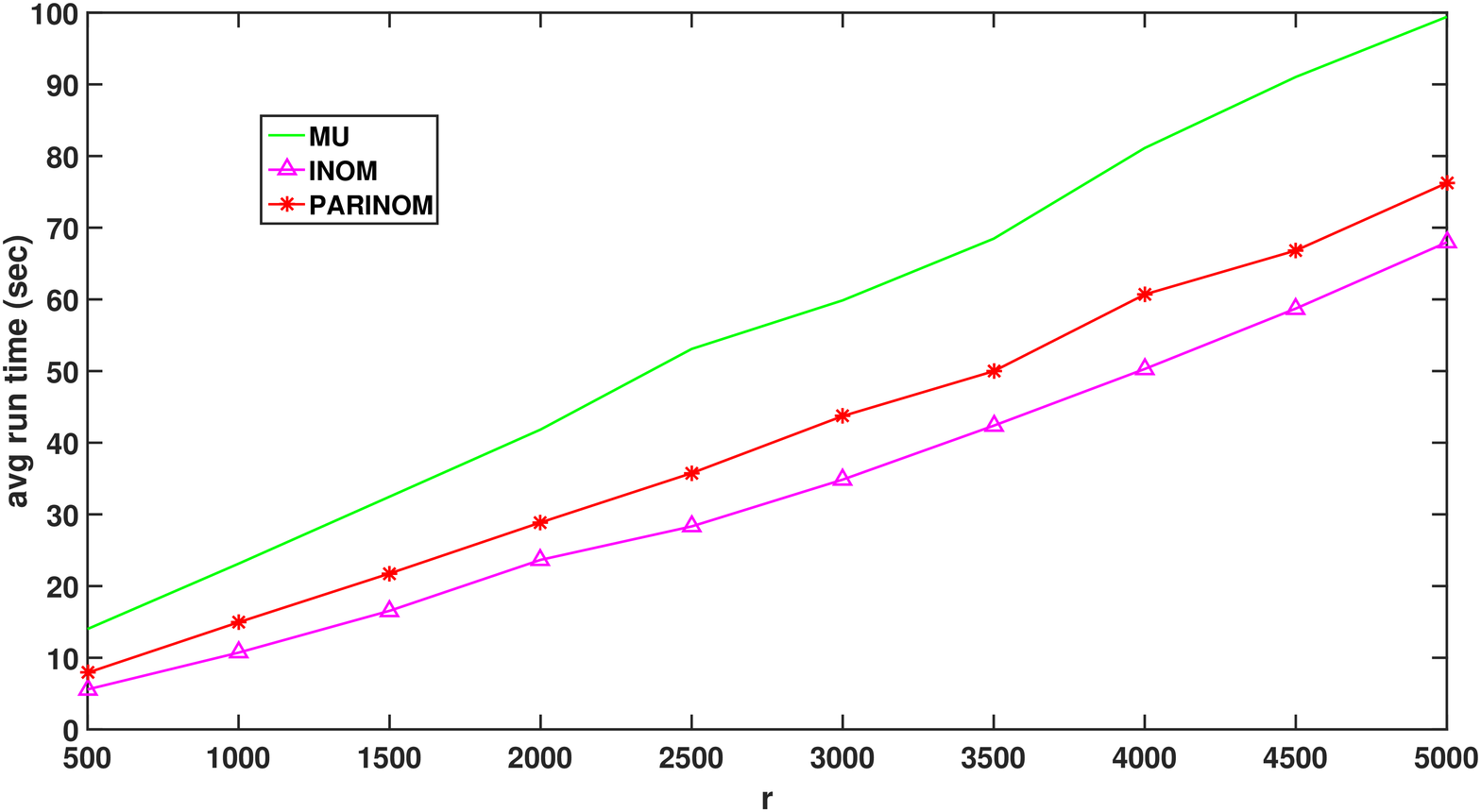}
\caption{Comparison of run time of MU, INOM and PARINOM}
\end{subfigure}
\caption{Comparison of proposed algorithm with existing algorithms by varying $r$, $m=50000$ and $n=10000$ and when $\bv$ is a $70 \%$  sparse matrix}
\label{sparse_matrix}
\end{figure}

\begin{figure}[H]
\centering
\begin{subfigure}{0.49\textwidth}
\centering
\captionsetup{justification=centering}
\includegraphics[height=2.5in,width=3.3in]{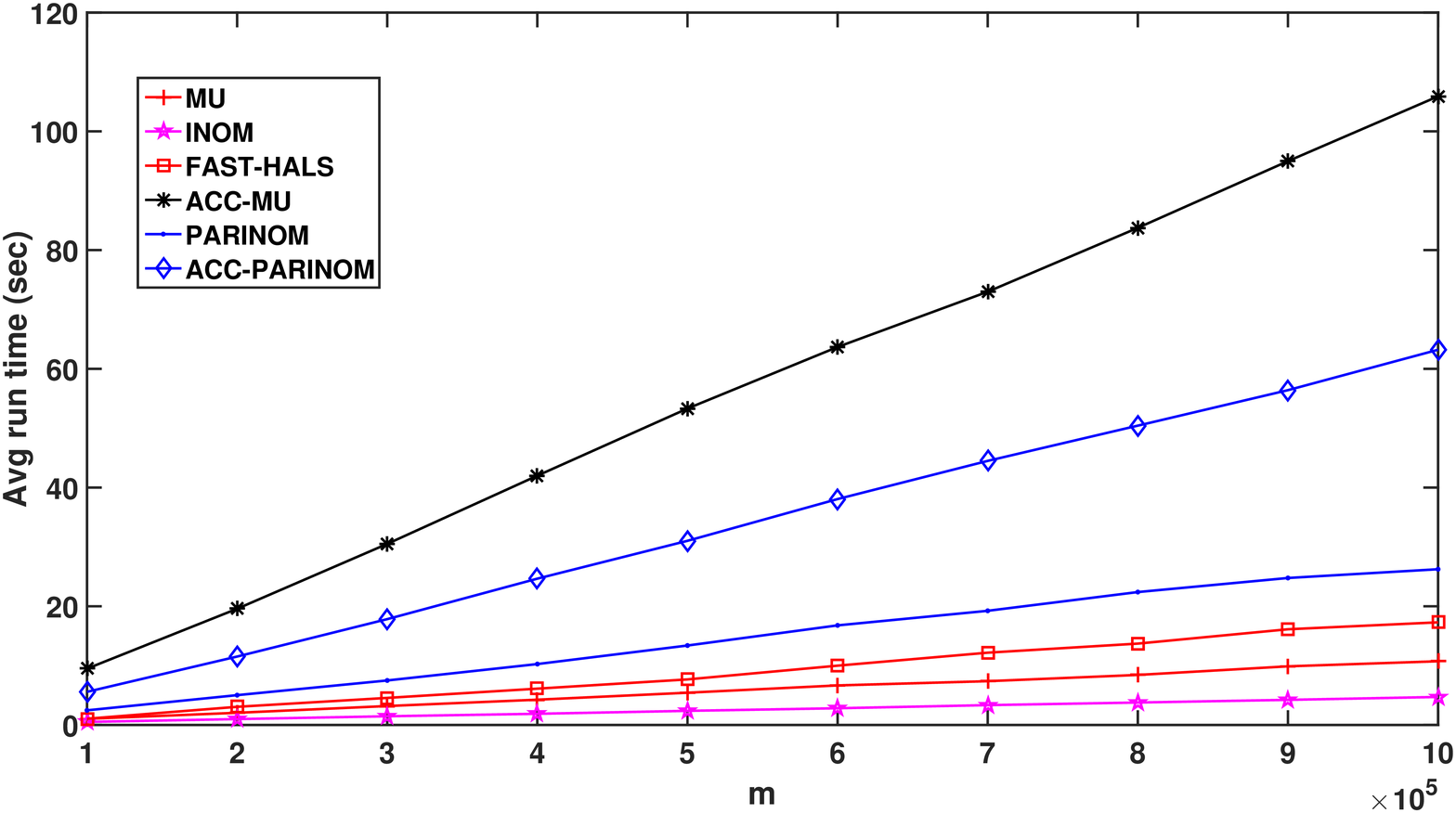}
\caption{Comparison of run time of proposed algorithms with Fast-Hals, MU and Acc-MU.}
\end{subfigure}
\begin{subfigure}{0.49\textwidth}
\centering
\captionsetup{justification=centering}
\includegraphics[height=2.5in,width=3.3in]{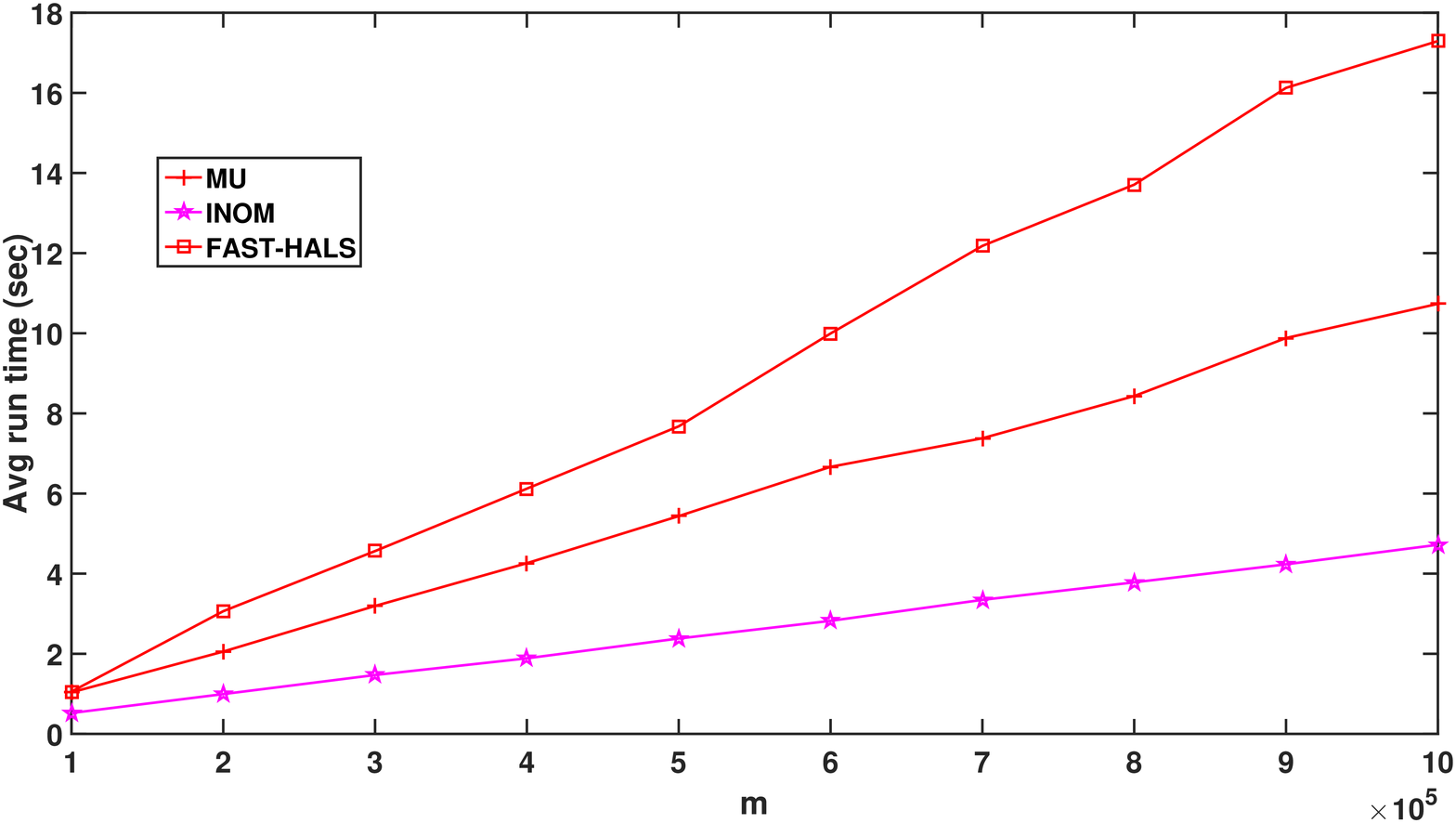}
\caption{Comparison of run time of MU, INOM and Fast-Hals}
\end{subfigure}
\caption{Comparison of proposed algorithm with existing algorithms for varying $r$, $m=50000$ and $n=10000$ when $\bv$ is dense matrix. }
\label{vary_m_dense_matrix}
\end{figure}

\begin{figure}[H]
\centering
\begin{subfigure}{0.49\textwidth}
\centering
\captionsetup{justification=centering}
\includegraphics[height=2.5in,width=3.3in]{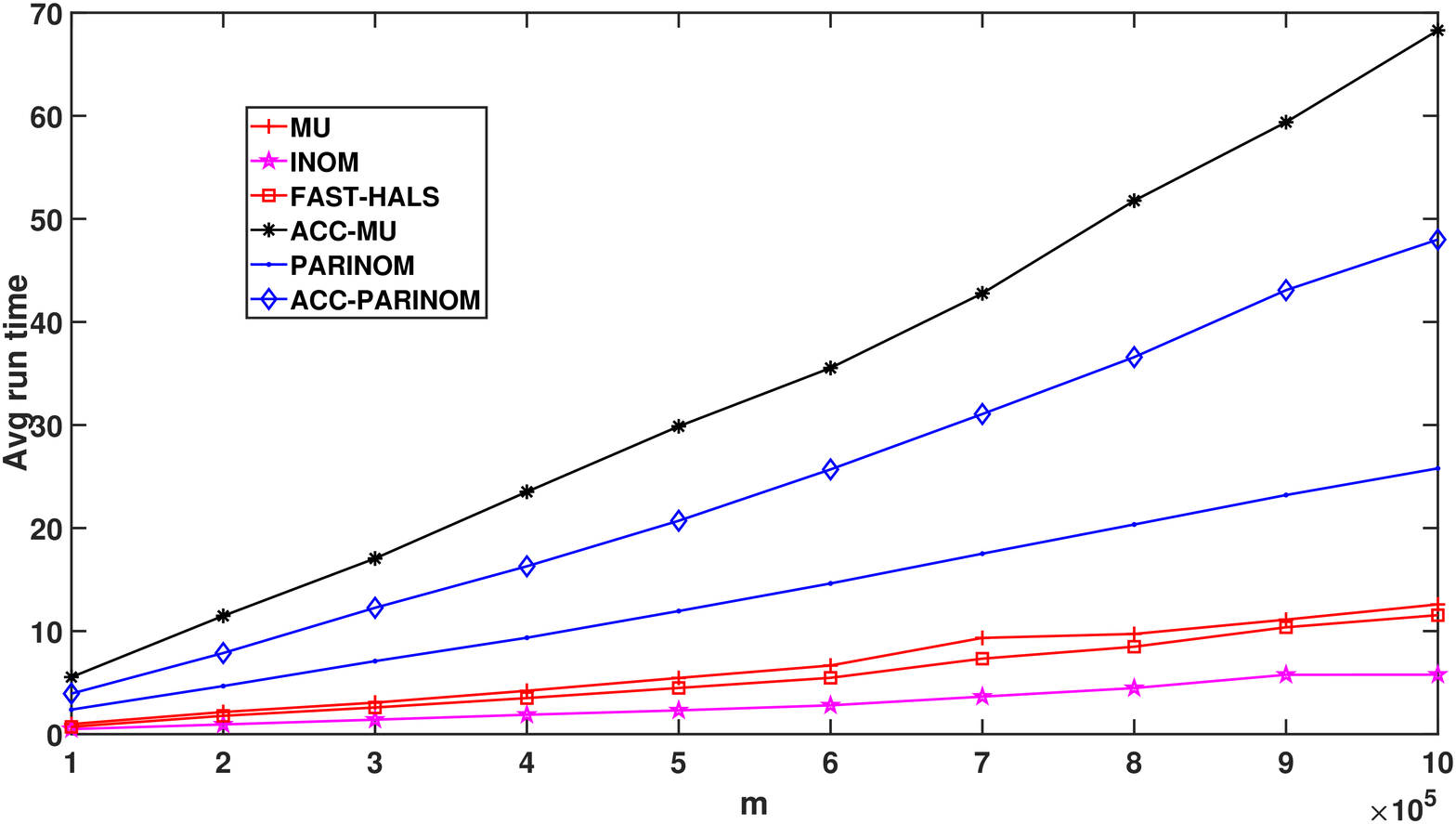}
\caption{Comparison of run time of proposed algorithms with Fast-Hals, MU and Acc-MU.}
\end{subfigure}
\begin{subfigure}{0.49\textwidth}
\centering
\captionsetup{justification=centering}
\includegraphics[height=2.5in,width=3.3in]{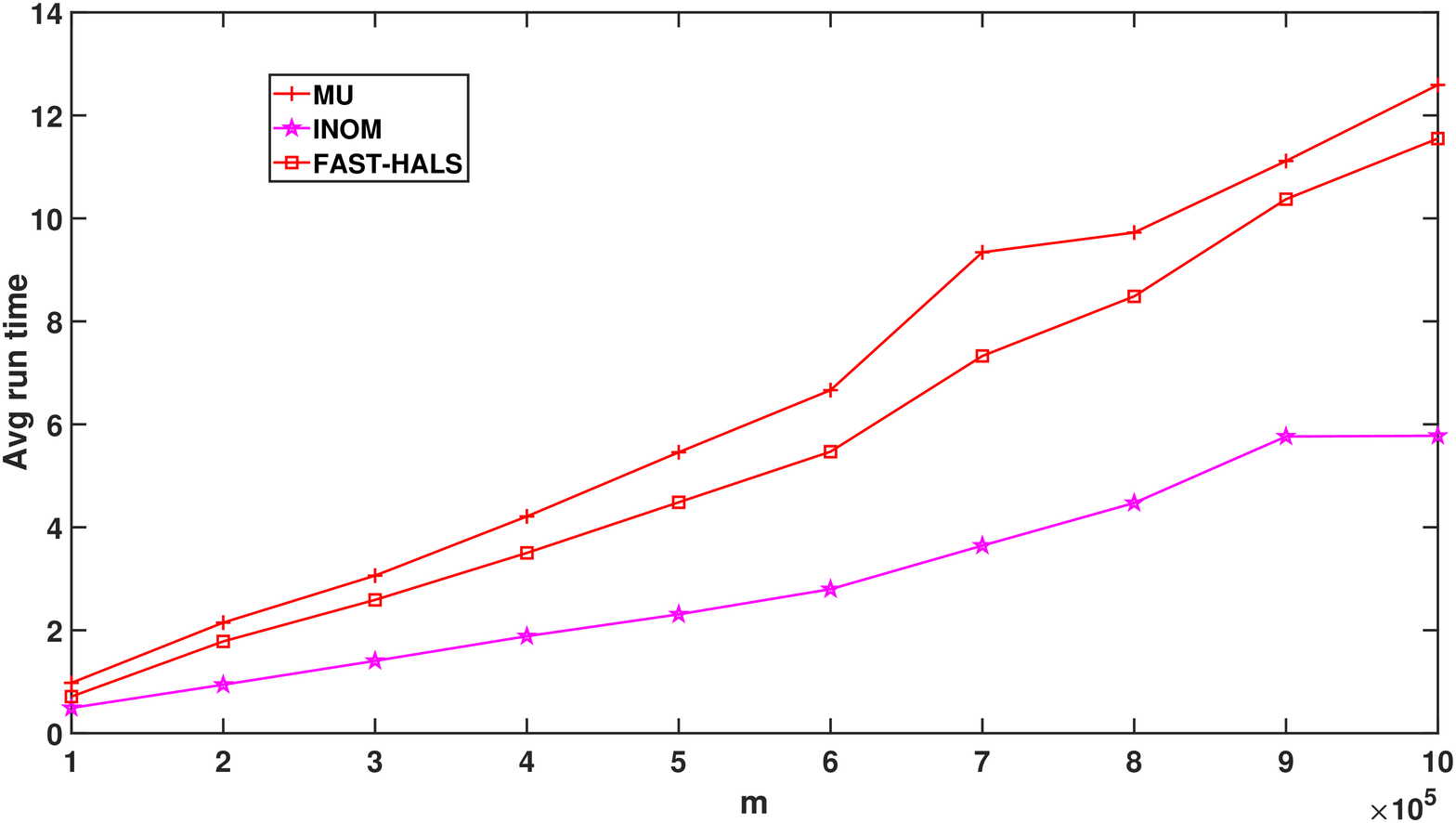}
\caption{Comparison of run time of MU, INOM and Fast-Hals}
\end{subfigure}
\caption{Comparison of proposed algorithm with existing algorithms for varying $r$, $m=50000$ and $n=10000$ when $\bv$ is $70 \%$ sparse matrix.}
\label{vary_m_sparse_matrix}
\end{figure}

From Fig. \ref{fig1}(a), it can be seen that the proposed algorithms are monotonic in nature. Accelerated \textbf{PARINOM} has faster convergence when compared to \textbf{PARINOM} algorithm. However, due to the additional steps involved in accelerated \textbf{PARINOM}, the per iteration cost is more in Accelerated \textbf{PARINOM} when compared to \textbf{PARINOM}. From Fig. \ref{fig1} (b), it can be seen that \textbf{INOM} algorithm takes lesser time to converge when compared to the rest of the algorithms. \\
\\
2.a. In this simulation, we compare the proposed algorithms with Fast-HALS, MU and accelerated MU for a dense matrix $\bv$ with $m=50000$ and $n =10000$. Dense matrix factorization has application in image processing and in video analysis. The elements of $\bv$ matrix was randomly generated from a uniform distribution from [100, 200]. Initial values of $\bw$ and $\bh$ was randomly generated from a uniform distribution from [0, 1]. The comparison is done based on how quickly the algorithms reduce the initial objective value $f_{_{\textrm{NMF}}}(\bw^{0},\bh^{0})$ to about $70\%$ of the initial objective value $f_{_{\textrm{NMF}}}(\bw^{0},\bh^{0})$ for different values of $r$- which was varied from $500$ to $5000$ in steps of $500$.  The columns of $\bv$ and $\bw$ were normalized. The initial objective value for all the algorithms were kept same. The run time was averaged over $50$ trials. Fig. \ref{large_matrix} (a) compares the proposed algorithms with Fast-HALS, MU and accelerated MU algorithm for $m=50000$ and $n=10000$.
%
%
%
%
From Fig. \ref{large_matrix} (a), it can be seen that as the size of $r$ increases, Fast-HALS takes the most time to reduce the initial  objective value $f_{_{\textrm{NMF}}}(\bw^{0},\bh^{0})$ to about $70\%$ of the initial objective value $f_{_{\textrm{NMF}}}(\bw^{0},\bh^{0})$. Fig.\ref{large_matrix} (b) compares MU, \textbf{INOM} and \textbf{PARINOM}. From this figure, it can be seen that \textbf{INOM} takes the least time.\\
\\
2.b. We repeat the above simulation for a sparse matrix $\bv$. Sparse matrix factorization has application in text mining. The elements of a $70 \%$ sparse matrix was randomly generated from a normal distribution with negative elements pushed to zero. Fig.\ref{sparse_matrix} (a) compares the proposed algorithms with Fast-Hals, MU and accelerated MU algorithm for $m=50000$ and $n=10000$ when the matrix $\bv$ is $70 \%$ sparse. 
From Fig.\ref{sparse_matrix} (a), it can be seen that \textbf{INOM} and \textbf{PARINOM} takes the least time when compared to the state-of-the art algorithms. Also, as the size of $r$ increases, Fast-HALS takes the most time to reduce the initial objective value $f_{_{\textrm{NMF}}}(\bw^{0},\bh^{0})$ to $70\%$ of the initial objective value $f_{_{\textrm{NMF}}}(\bw^{0},\bh^{0})$. Fig.\ref{sparse_matrix} (b) shows only the performance of \textbf{INOM}, \textbf{PARINOM} and MU for better readability. \\
\\
3.a. In this simulation, we vary the size of $\bv$ and compare the performance of the algorithm. The comparison is done based on how quickly the algorithms reduce the initial objective value $f_{_{\textrm{NMF}}}(\bw^{0},\bh^{0})$ to about $70\%$ of the initial objective value $f_{_{\textrm{NMF}}}(\bw^{0},\bh^{0})$. $m$ was varied from $100000 $ to $1000000$ in steps of $100000$ and $n$ and $r$ were equal to $1000$ and $100$, respectively. The elements of $\bv$ matrix was randomly generated from a uniform distribution from [100, 200]. Initial values of $\bw$ and $\bh$ was randomly generated from a uniform distribution from [0, 1]. The columns of $\bv$ and $\bw$ were normalized. The initial objective value for all the algorithms were kept same. The run time was averaged over $50$ trials. Fig. \ref{vary_m_dense_matrix} (a) compares the proposed algorithms with Fast-Hals, MU and accelerated MU algorithm for $r=100$ and $n=1000$. Fig. \ref{vary_m_dense_matrix} (b) compares the performance of \textbf{INOM}, MU and Fast-Hals.
From Fig.\ref{vary_m_dense_matrix} (b), it can be seen that \textbf{INOM} takes the least time when compared to the state-of-the art algorithms. \\
\\
3.b. The above simulation is repeated for a sparse matrix $\bv$. The elements of a $70 \%$ sparse matrix was randomly generated from a normal distribution with negative elements pushed to zero. Fig.\ref{vary_m_sparse_matrix} (a) compares the proposed algorithms with Fast-Hals, MU and accelerated MU algorithm for $r=100$ and $n=1000$. Fig.\ref{vary_m_sparse_matrix} (b) compares the performance of \textbf{INOM}, MU and Fast-Hals. From Fig.\ref{vary_m_sparse_matrix} (b), it can be seen that \textbf{INOM} performs better than MU and Fast-Hals.  

\begin{figure}[H]
\centering
\begin{subfigure}{0.49\textwidth}
\centering
\includegraphics[height=2.5in,width=3.3in]{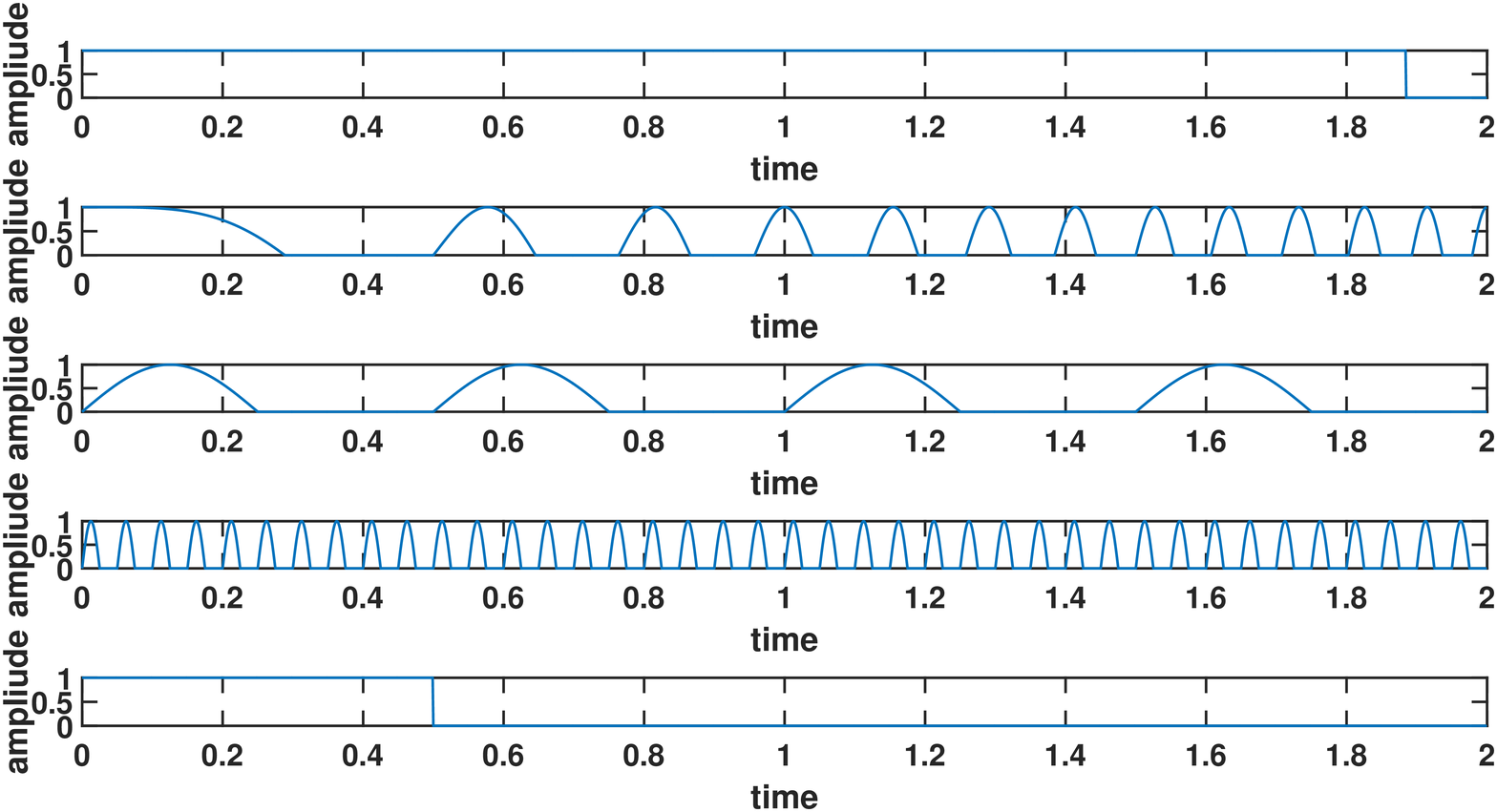}
\caption{}
\end{subfigure}
\begin{subfigure}{0.49\textwidth}
\centering
\includegraphics[height=2.5in,width=3.3in]{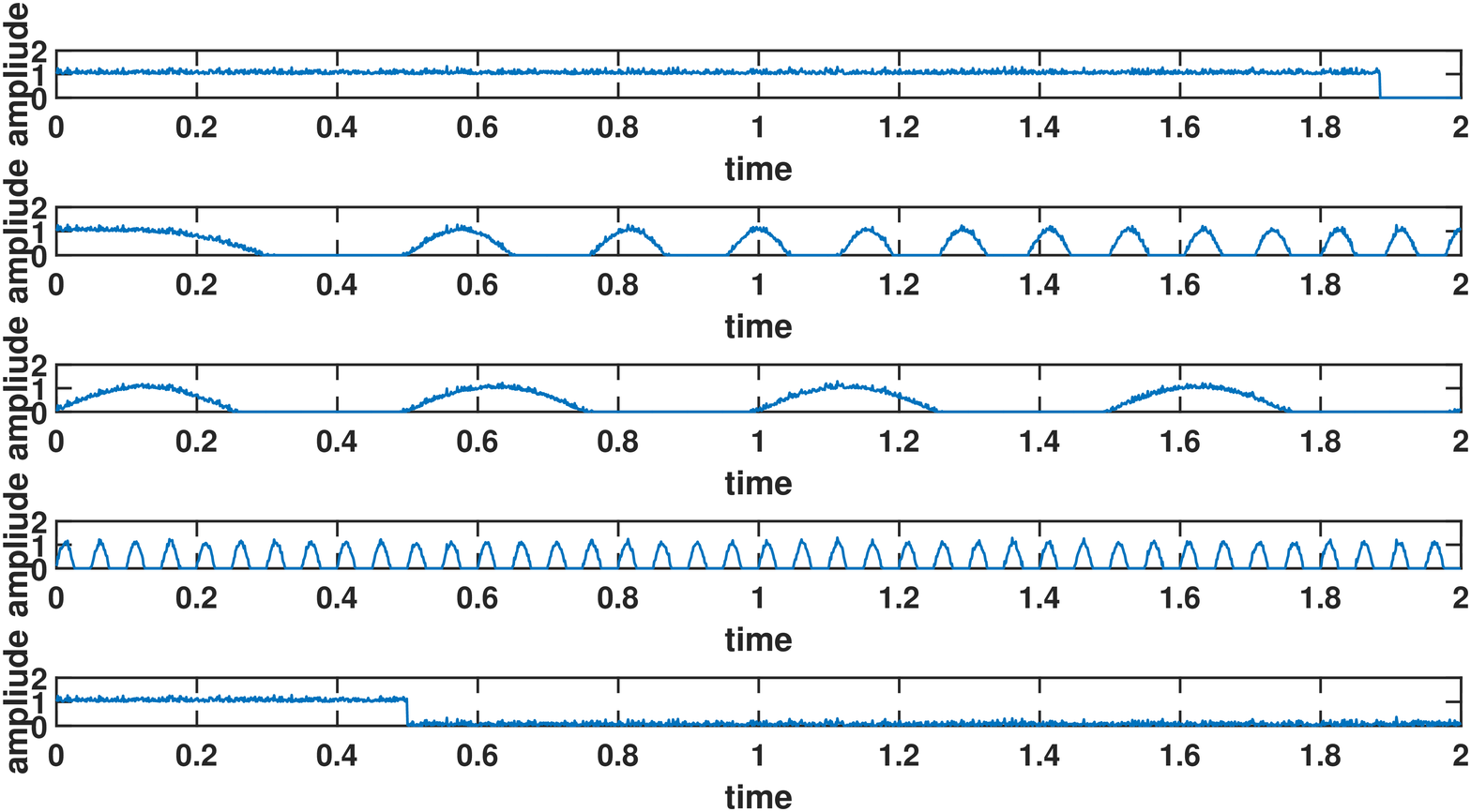}
\caption{}
\end{subfigure}
\caption{(a) Source signal without noise (b) Noisy source signal with variance = $0.1$}
\label{fig4}
\end{figure}
\begin{figure}[H]
\centering
\begin{subfigure}{0.49\textwidth}
\centering
\includegraphics[height=2.5in,width=3.3in]{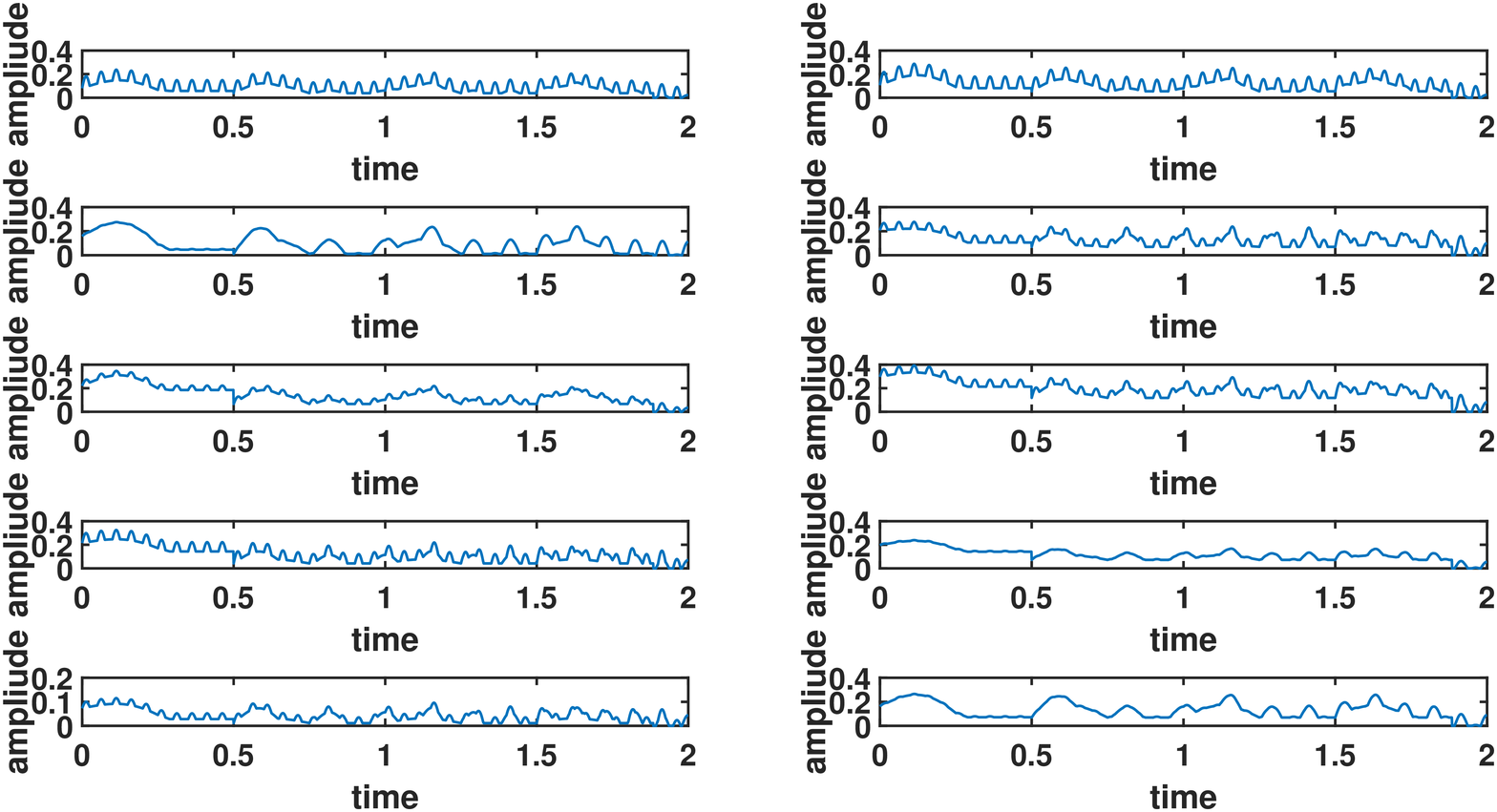}
\caption{}
\end{subfigure}
\begin{subfigure}{0.49\textwidth}
\centering
\includegraphics[height=2.5in,width=3.3in]{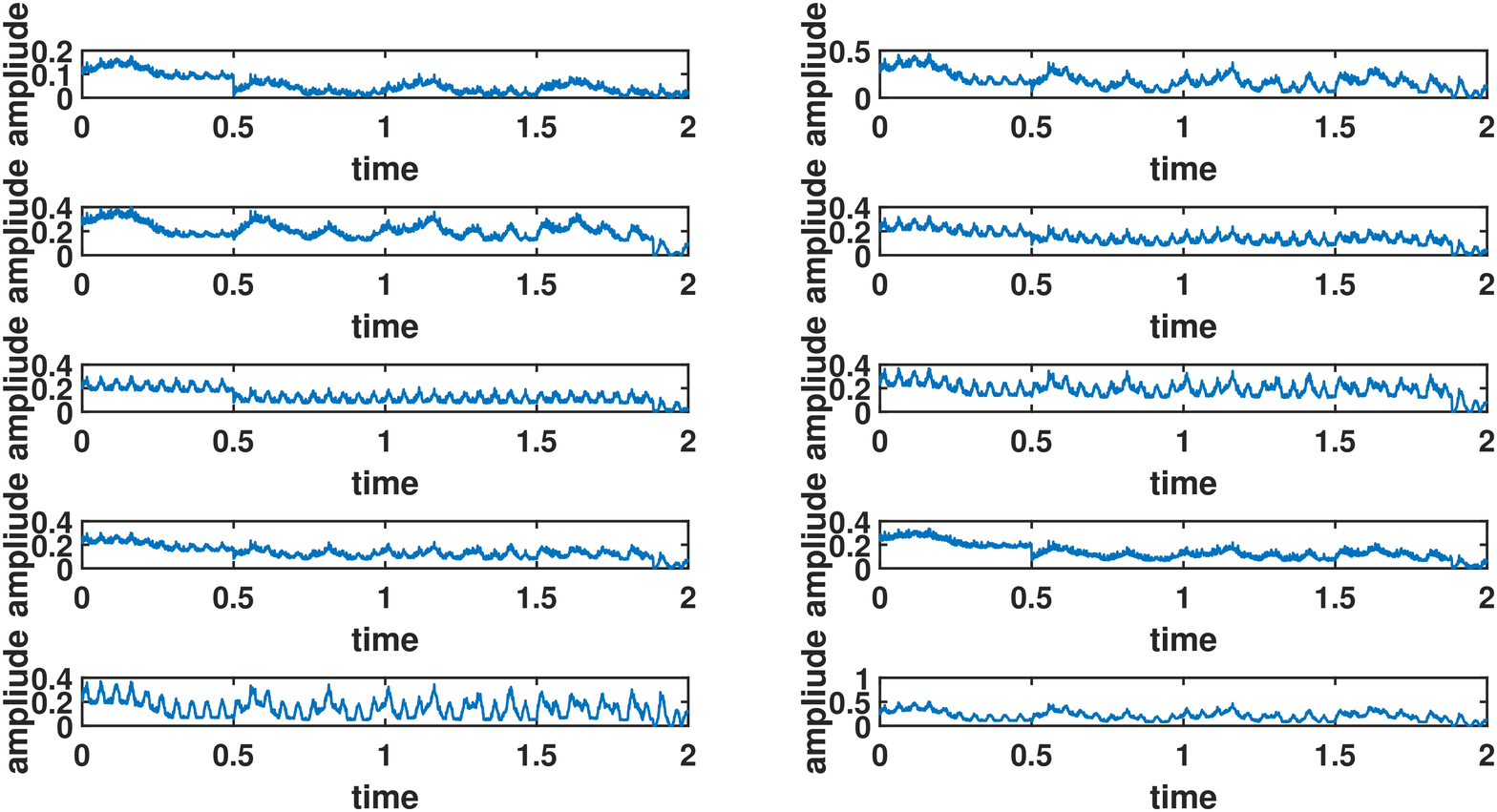}
\caption{}
\end{subfigure}
\caption{(a) Observed signal when source signal has no noise. (b) Observed signal when source signal has noise variance = $0.1$}
\label{fig5}
\end{figure}

4. In this simulation, an application of Nonnegative Matrix Factorization on Blind Source Separation(BSS) is shown. In BSS, source signals have to be separated from a set of mixed signals without knowing the mixing process. It is commonly used in audio signal processing, image processing, biomedical signal processing and in digital communications. In the latter case, the $\bw$ matrix can be thought as the channel response and each row of $\bh$ consists of the source signal send by the transmitter array. The receiver sensor array then receives a linear mixture of the source signal. The task of BSS is to reconstruct the source signal from the received signal.  \\
\\
Five source signals for $10$ seconds were simulated-a square wave, a rectangular wave, two sine waves of frequency $2$ Hz and $20$ Hz and a chirp signal which begins with $0$ Hz at t = $0$ sec and cross $30$ Hz after $t=5$ sec. The negative values of the source signals were made zero. $\bw$ matrix of size $200$ $\times$ $5$ was randomly generated. To evaluate the proposed algorithm in presence of noise, random noise with variance $0.01$ was added to the five source signals. Source signals and some of the observed signals are shown in Fig. \ref{fig4} and Fig. \ref{fig5} respectively.  

\begin{figure}[H]
\centering
\begin{subfigure}{0.49\textwidth}
\centering
\includegraphics[height=2.5in,width=3.3in]{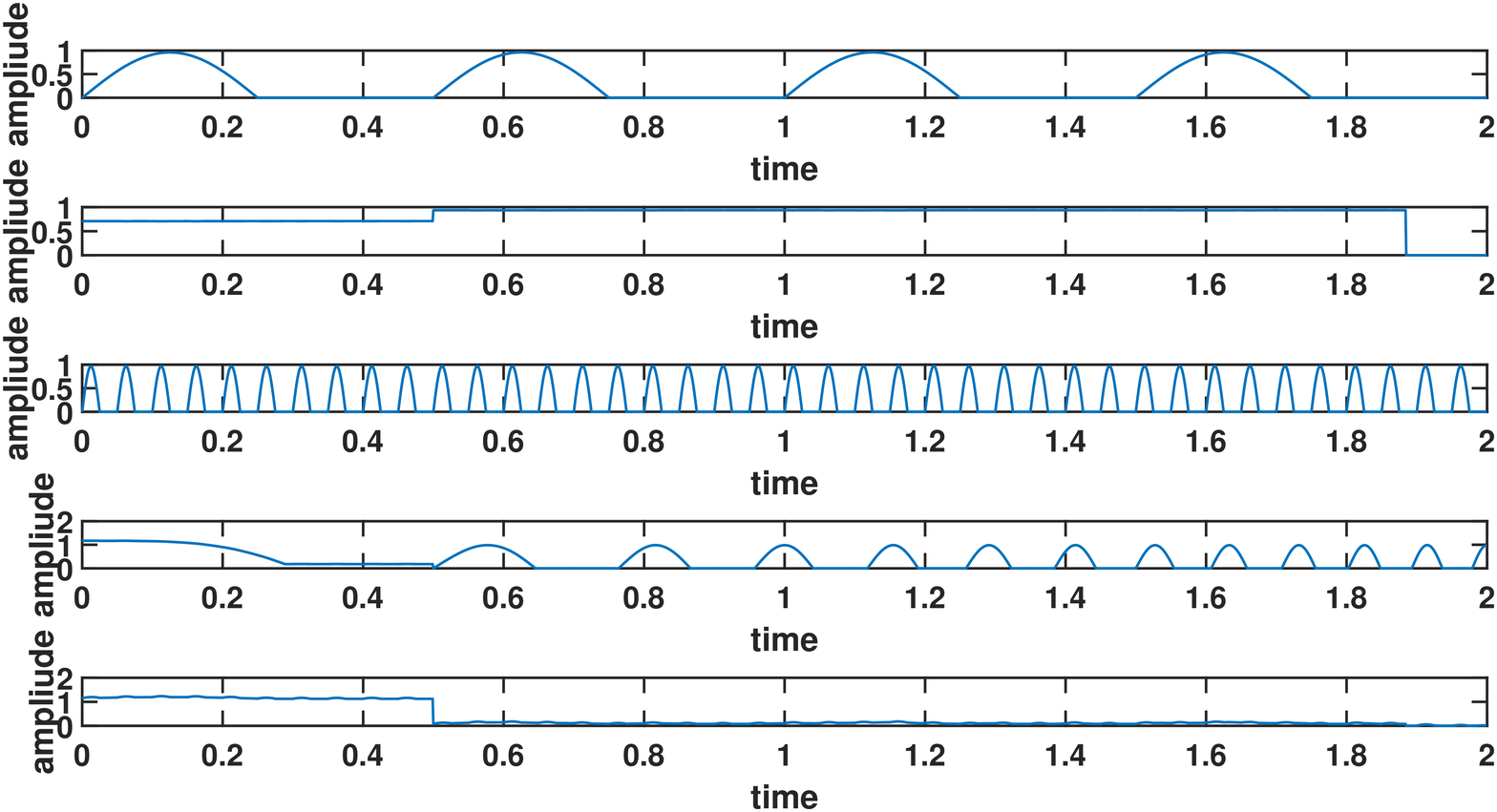}
\caption{}
\end{subfigure}
\begin{subfigure}{0.49\textwidth}
\centering
\includegraphics[height=2.5in,width=3.3in]{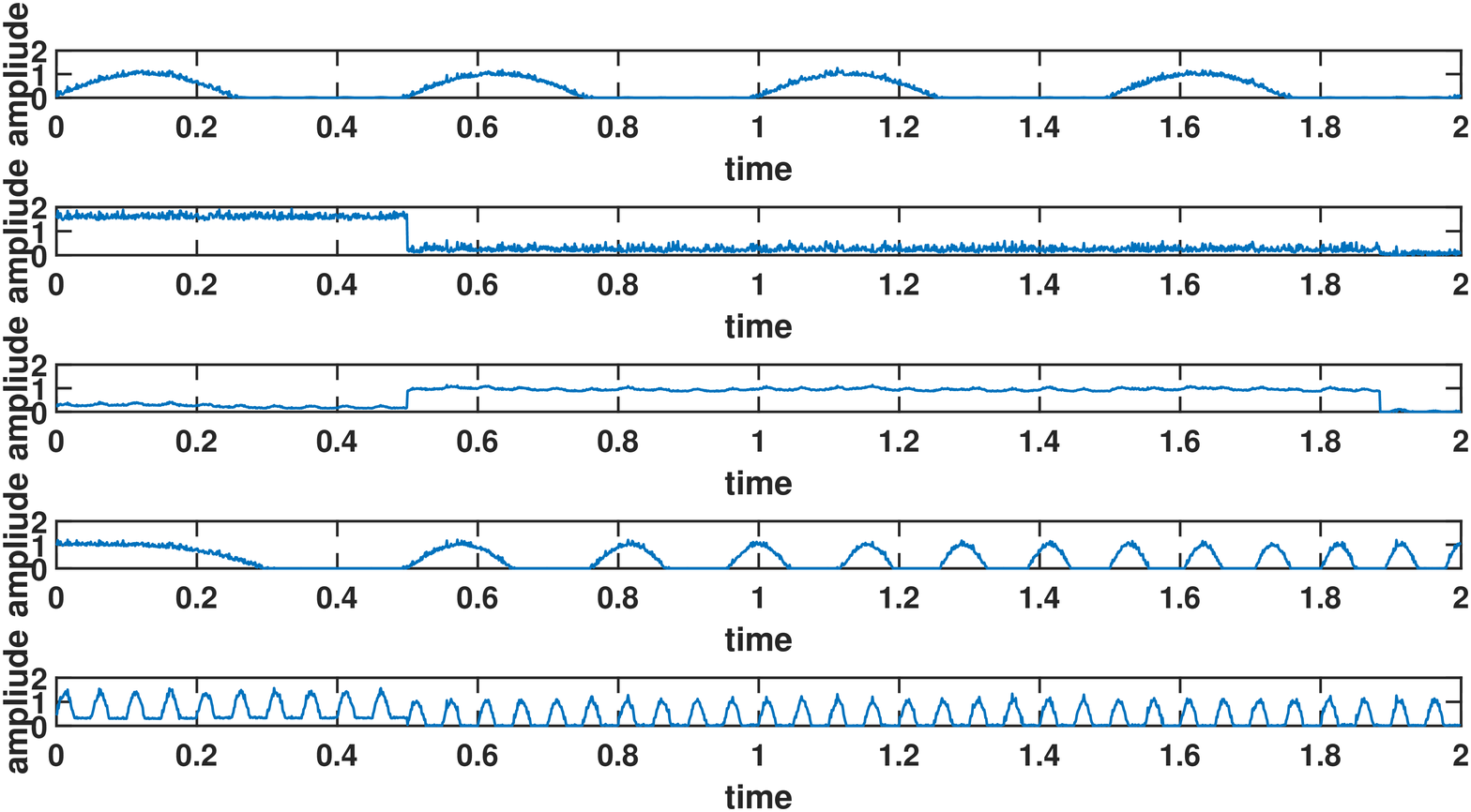}
\caption{}
\end{subfigure}
\caption{Reconstructed signal by \textbf{INOM} (a) without noise. (b) with noise variance = $0.1$}
\label{fig6}
\end{figure}

\begin{figure}[H]
\centering
\begin{subfigure}{0.49\textwidth}
\centering
\includegraphics[height=2.5in,width=3.3in]{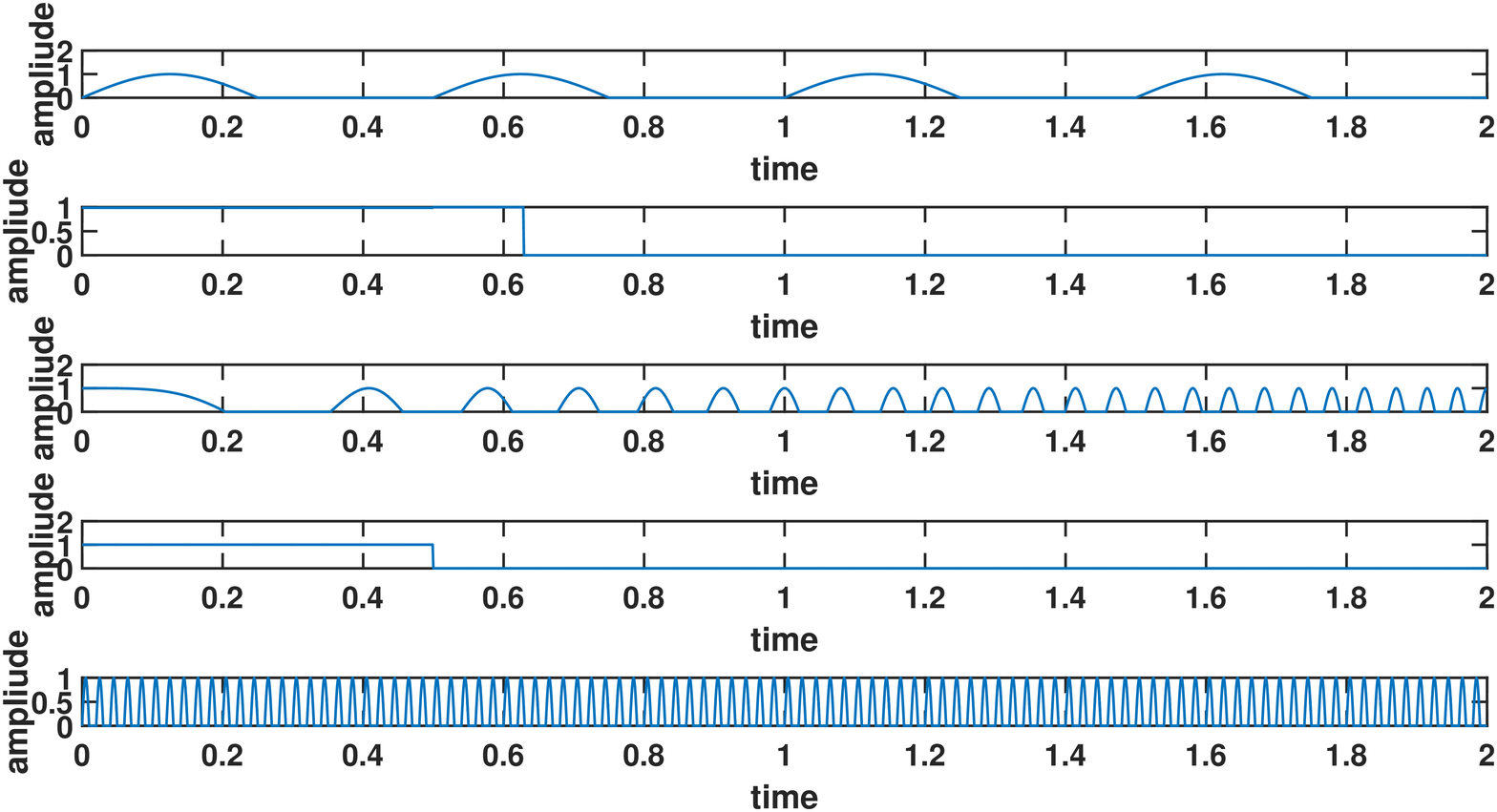}
\caption{}
\end{subfigure}
\begin{subfigure}{0.49\textwidth}
\centering
\includegraphics[height=2.5in,width=3.3in]{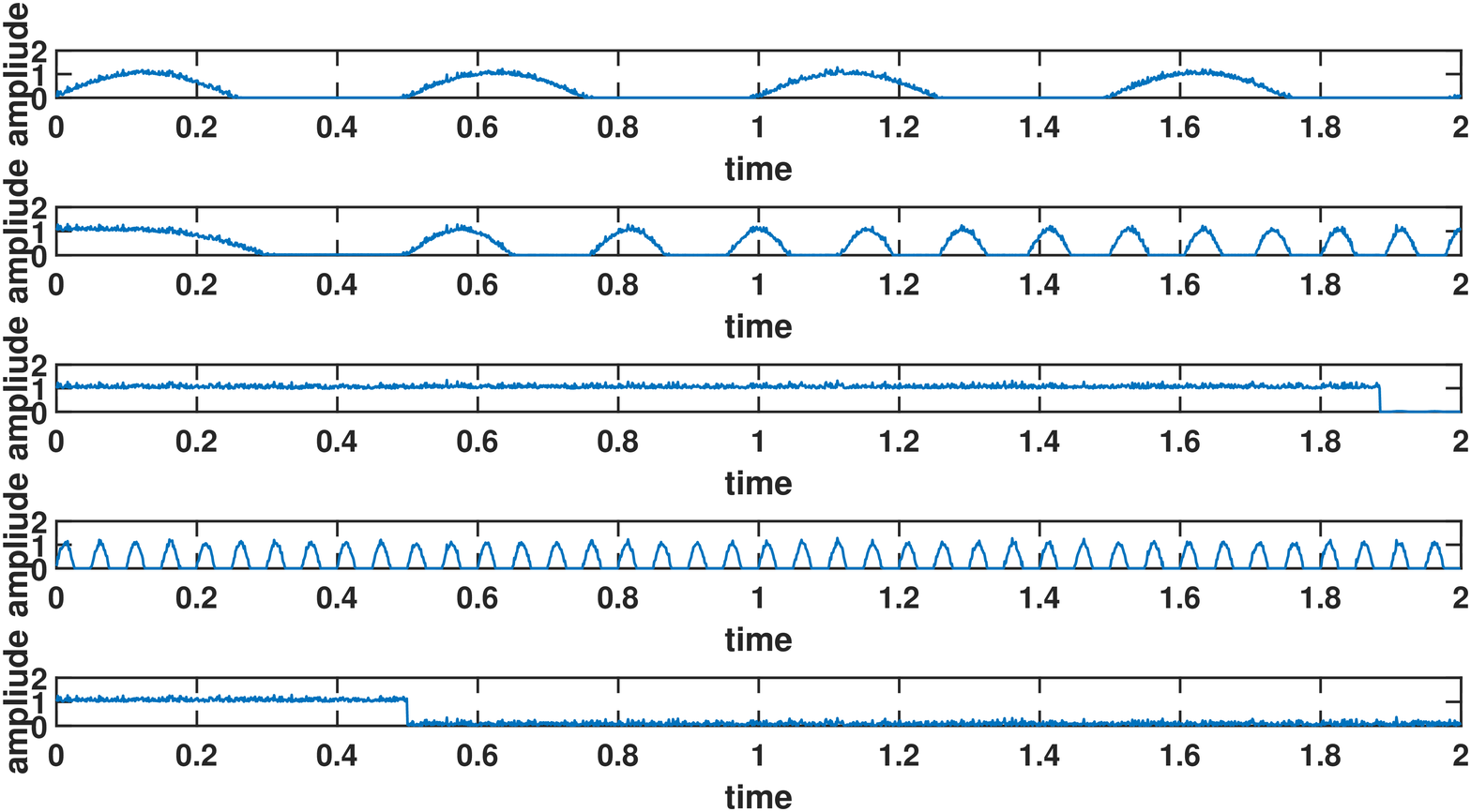}
\caption{}
\end{subfigure}
\caption{Reconstructed signal by Fast-Hals (a) without noise (b) with noise variance = $0.1$}
\label{fig7}
\end{figure}

\begin{figure}[H]
\centering
\begin{subfigure}{0.49\textwidth}
\centering
\includegraphics[height=2.5in,width=3.3in]{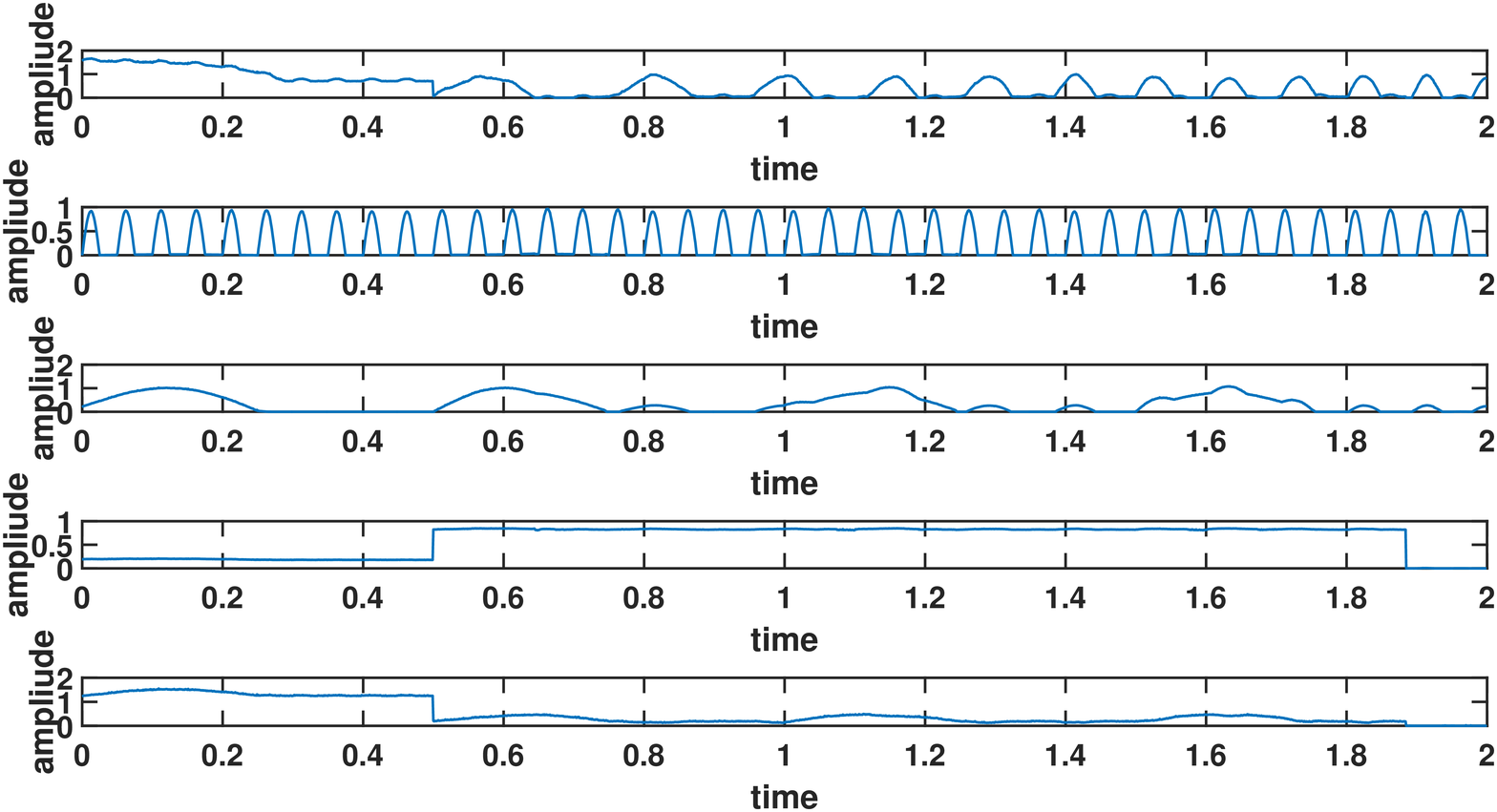}
\caption{}
\end{subfigure}
\begin{subfigure}{0.49\textwidth}
\centering
\includegraphics[height=2.5in,width=3.3in]{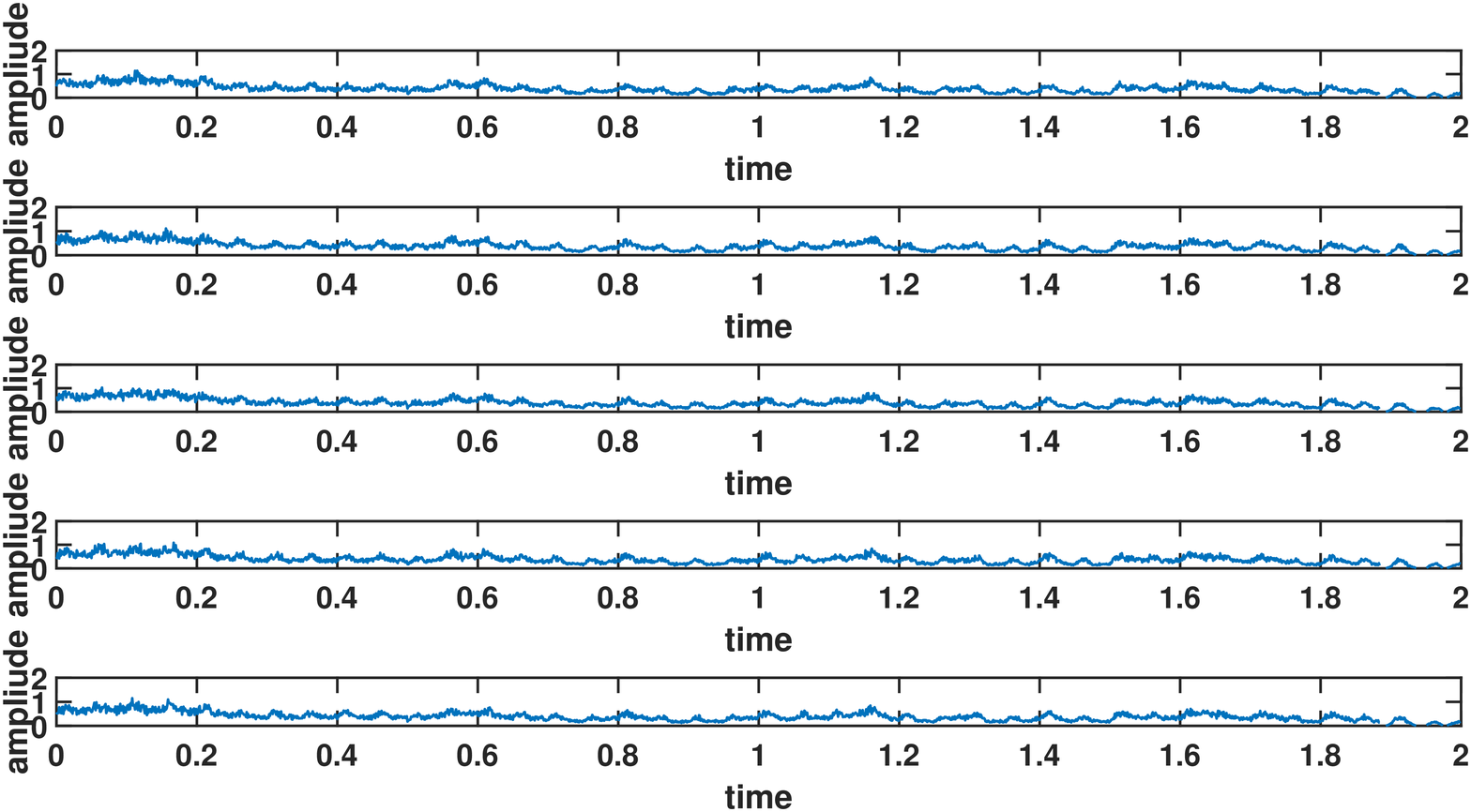}
\caption{}
\end{subfigure}
\caption{Reconstructed signal by MU (a) without noise (b) with noise variance = $0.1$}
\label{fig8}
\end{figure}

$\bw$ and $\bh$ were initialized randomly from Uniform distribution from [100, 500] and [200, 400], respectively. The algorithms were made to run till 1000 iterations or unless the relative change in cost function was $10^{-8}$. Figure. \ref{fig6} and \ref{fig7} shows the reconstructed signal by \textbf{INOM} and Fast-Hals, respectively. MU algorithm was not able to reconstruct the signal in the presence of noise and the same is shown in Fig. \ref{fig8}.

\section{Possible Extension and Conclusion}\label{sec:7}
In this paper, we proposed two different MM based algorithms - \textbf{INOM} and \textbf{PARINOM} to solve the NMF  problem. \textbf{INOM} sequentially updates the $\bw$ and $\bh$ matrices, while \textbf{PARINOM} parallely updates them. The update equations of $\bw$ and $\bh$ matrices in case of \textbf{INOM} resembles the update equation of  Gradient Descent algorithms with adaptive step sizes. We also prove that the proposed algorithms are monotonic and converge to the stationary point of the NMF problem. Various computer simulations were performed to compare the proposed algorithms with the existing algorithms. It was found that \textbf{INOM} performed better than existing algorithms. \\
\\
The effect of parallel update will be predominantly visible when one has to decompose a multidimensional array into matrices of varied dimensions. This problem is called NonNegative Tensor Factorization (NTF) \cite{ntf_recent_1}, \cite{ntf_recent_3}, \cite{ntf_recent_2}, \cite{sidiropoulos_1}. A brief overview of NTF is given below.\\
\\
Suppose, a tensor {\bf{X}} of order N $\mathbb{R}^{{m_{1} \times {m_{2}}}\cdots {m_{N}}}$ is given. This can be represented as sum of rank one tensors i.e. 
\begin{equation}
\begin{array}{ll}
\mathbf{X} = \displaystyle\sum_{k=1}^{K} {\mathbf{a_{k}}}^{(1)} \otimes {\mathbf{a_{k}}}^{(2)}\cdots {\mathbf{a_{k}}^{(N)}}
\end{array}
\end{equation}
where ${\mathbf{a_{k}}}^{(1)}$ belongs to $\mathbb{R}^{m_{1}}$, ${\mathbf{a_{(k)}}}^{(2)}$ belongs to $\mathbb{R}^{m_{2}}$ and so on. K is rank of the tensor and $\otimes$ represents outer product. This can be written in a compact form as:
\begin{equation}
\begin{array}{ll}
\mathbf{X} = [\mathbf{A}^{(1)}, \mathbf{A}^{(2)} \cdots \mathbf{A}^{(N)}]
\end{array}
\end{equation}
where 
${\mathbf{A}}^{(1)} = [{\mathbf{a_{1}}}^{(1)}, {\mathbf{a_{2}}}^{(1)},\cdots,{\mathbf{a_{K}}}^{(1)}]$ is of size $m_{1} \times K$. Similarly, ${\mathbf{A}}^{(2)} = [{\mathbf{a_{1}}}^{(2)}, {\mathbf{a_{2}}}^{(2)},\cdots,{\mathbf{a_{K}}}^{(2)}]$ and is of size $m_{2} \times K$. $[\mathbf{A}^{(1)}, \mathbf{A}^{(2)} \cdots \mathbf{A}^{(N)}]$ represents $\sum_{k=1}^{K} {\mathbf{a_{k}}}^{(1)} \otimes {\mathbf{a_{k}}}^{(2)}\cdots {\mathbf{a_{k}}^{(N)}}$. Hence, the NTF problem can be formulated as:
\begin{equation} \label{NTF}
\begin{array}{ll}
\textrm{NTF:} \quad \underset{\mathbf{A}^{(1)},\mathbf{A}^{(2)}\cdots\mathbf{A}^{(N)} \geq 0}{\rm minimize} \: \left(\|\mathbf{X} - [\mathbf{A}^{(1)}, \mathbf{A}^{(2)} \cdots \mathbf{A}^{(N)}]\|_F^2\right) 
\end{array}
\end{equation}
When compared to parallel update of matrices, alternating minimization will definitely take significantly more time to solve the problem. Since NTF is a generalization of NMF, the proposed algorithms can also be applied to solve NTF and the parallel update algorithm - \textbf{PARINOM} would be more handy for NTF. \\
\\
\bibliographystyle{IEEEtran} 
\bibliography{ref_nnmf}

\end{document}